\documentclass[preprint,12pt]{elsarticle}

\usepackage{amssymb}
\usepackage{amsmath}
\usepackage{epsfig}
\usepackage{hyperref}
\usepackage{verbatim}
\allowdisplaybreaks

\newtheorem{Proposition}{Proposition}[section]
  \newtheorem{Remark}[Proposition]{Remark}
  \newtheorem{Corollary}[Proposition]{Corollary}
  \newtheorem{Lemma}[Proposition]{Lemma}
  
  \newtheorem{Theorem}{Theorem}[section]
 
 \newtheorem{Definition}[Proposition]{Definition}

 \newproof{pf}{Proof}

\journal{Journal of Differential Equations}

\begin{document}

\begin{frontmatter}



\title{Existence, Uniqueness, Analyticity, and Borel Summability for 
Boussinesq Equations}


\author[ad1]{H. Rosenblatt}
\author[ad2]{S. Tanveer}

\address[ad1]{Department of Mathematics: The Ohio State University \\
              529 Math Tower 231 West 18th Avenue Columbus, OH 43210-1174\\
              email: rosenblatt@math.ohio-state.edu phone:1-614-292-1923}
              
\address[ad2]{402 Math Tower 231 West 18th Avenue Columbus, OH 43210-1174\\
               email: tanveer@math.ohio-state.edu}

\begin{abstract}
Through Borel summation methods, 
we analyze the 
Boussinesq
equations for coupled fluid velocity and temperature fields:
\begin{align}\label{star}
u_t-\nu\Delta u&=-P[ u\cdot \nabla u - a e_2\Theta] +f\\ \nonumber
\Theta_t-\mu \Delta \Theta&= -u\cdot \nabla \Theta.
\end{align}
We prove that 
an equivalent system of integral equations in the Borel variable $p\in \mathbb{R}^{+}$ dual to $1/t$ has a unique solution in a class of 
exponentially bounded functions, implying
the existence of a classical solution to (\ref{star}) in
a complex $t$-region that includes a real positive time axis segment.
For analytic initial data and forcing, it is shown that the solution is 
Borel summable, implying that
that formal series in powers of $t$ is Gevrey-1 asymptotic, and within the time interval of existence,
the solution remains analytic with the same analyticity strip width as the initial data and forcing. We also determine conditions
on the integral equation solution that improve the estimate for existence time.

\end{abstract}


\end{frontmatter}


\section{Introduction}
\label{intro}
We consider the Boussinesq equations for coupled fluid velocity and temperature fields
derived under the assumption that the temperature induced density has
negligible effect on momentum but causes  
a significant buoyant force. The corresponding evolution equations
for $u:\mathbb{R}^d\times\mathbb{R}^+\rightarrow \mathbb{R}^d$ 
and $\Theta:\mathbb{R}^d\times\mathbb{R}^+\rightarrow \mathbb{R}$ for dimension $d=2, 3$ 
in non-dimensional form are:     
\begin{align}
\label{B}
u_t-\nu\Delta u&=-P[ u\cdot \nabla u - a e_2\Theta] +f ~~,~~~~~u(x,0) = u_0 (x)
\\ \nonumber
\Theta_t-\mu \Delta \Theta&= -u\cdot \nabla \Theta ~~~,~~~\Theta (x, 0) =\Theta_0 (x) 
\end{align}
where $P=I-\nabla \Delta^{-1}(\nabla\cdot )$ is the Hodge projection operator to 
the space of divergence free vector fields, $e_2$ is the unit vector 
aligned opposite to gravity, the parameter $a$ is proportional to gravity,
 and $(u, \Theta)$ are the nondimensional fluid velocity and temperature fields. 
We assume the initial conditions $u_0$ and the
forcing $f$ are divergence free and, for the sake of simplicity, assume $f$ to be time independent, 
although time dependence with some restrictions can be accommodated in a similar framework.
Using standard energy methods, see for instance \cite{Temam}, 
existence of Leray type solutions in $L^{\infty}(0,T,L^2(\mathbb{R}^d))\cap L^2(0,T,H^1(\mathbb{R}^d))$ 
follows easily 
for any $T>0$. In $\mathbb{R}^2$ a unique classical global 
solution can be shown to exist for all time 
In \cite{Cannon}, 
local existence and uniqueness for Boussinesq equation 
are shown in $L^p(0,T,L^q(\mathbb{R}^d))$ for $d<p<\infty$ and $\frac{d}{p}+\frac{2}{q}\leq 1$. In $\mathbb{R}^3$ there is a unique solution 
under the additional assumption 
that the solution lies in $L^{\infty}(0,T,H^1(\mathbb{R}^3))$, see \cite{Cannon}. The case where $\mu=0$ has also been considered in the literature, and global well-posedness is proved proved in \cite{Boussinesq} for $2-d$.   

In the problem above, the existence of classical solutions, globally in time, remains an open problem 
as it is for the limiting ($a \rightarrow 0$) Navier
Stokes equation (NSE) in 3-D. 
Control of a higher order energy norm (like the $H^1$ norm of velocity) has remained 
a serious impediment despite extensive study of NSE. This
motivates one to look for alternate formulations of existence that do not rely on energy bounds.

The primary purpose of this paper is to show that the 
Borel based methods, developed earlier in \cite{smalltime} and \cite{longtime} in the context of  Navier-Stokes equation, 
can be extended to other evolutionary PDEs (partial differential equations) such as the Boussinesq equation. 
This provides an alternate existence and uniqueness
theory for a class of nonlinear PDEs. In this formulation, the question of global existence
of solution to the PDE becomes one of asymptotics 
for known solution to the associated nonlinear integral equations. 
While the asymptotics are still difficult, it is interesting to note
that an accelerated 
representation \cite{longtime} (see (\ref{accelerated}) in the ensuing)
for the related NSE results in a positive limiting kernel as 
$n \rightarrow \infty$, where  
majorization may be possible in terms of solution to a 
simpler integral equation.
We also show (Thm \ref{improved existence}) here 
how information about solution to the integral equation on a finite interval 
in the dual variable for specific initial condition and forcing may be used 
to obtain better exponential bounds in the Borel plane implying a longer existence time 
for classical solutions to the associated PDEs.

Borel summability has been an active area of research.
A vast literature has emerged recently in Borel summability theory, starting with the fundamental
contributions of Ecalle (see e.g. \cite{Ecalle} and \cite{Ecalle1}) whose consequences are far from being
fully explored, and it is impossible to give a quick account of the breadth of this field (See
for example \cite{Costin-Costin} for more references). There has also been
work in characterizing all small solutions for a generic system of ODEs \cite{Duke} or
difference equations \cite{Braaksma}. There has been work on  
PDEs as well, starting with linear equations \cite{Lutzetal}, \cite{Balser} 
followed by general results for  
a class of nonlinear system of PDEs in complex sectors \cite{CT01}, \cite{CT02}. A Borel based
approach has also led to
analysis of complex singularities for a specific PDE \cite{CT03}. 
Recent developments include Navier-Stokes initial value problem
(see \cite{longtime}, \cite{scripta}, \cite{CT1}). 
Recently \cite{Raza}, numerical schemes have been suggested 
for nonlinear PDEs, based on a Borel plane reformulation. Thus,
it is clear that the Borel based approach of the present paper is likely
to have both theoretical and practical value.
A bi-product of the present Borel based approach is that many analyticity properties of the PDE solution readily follow 
without additional arguments. For instance, the
time analyticity for $\Re \frac{1}{t}>\alpha$ follows
from (\ref{Laplacetransform1}) after noting the solution to the integral equation
is exponentially bounded in $p$.
While such analyticity results may also be obtained through other methods, 
see \cite{t-analyticity1} and \cite{t-analyticity2}, it follows more readily 
from the current method.
We also prove that 
the classical $H^2(\mathbb{R}^d)$ solution, which is unique, has 
the Laplace transform representation given here, 
provided initial data and forcing in the Fourier-space are in $L^1\cap L^{\infty}$. 
Furthermore, for analytic 
initial data and forcing, we prove that the formal expansion in powers of
$t$ is Borel summable and hence Gevrey-1 asymptotic for small $t$. 
As far as we know, these results are new and have not been obtained earlier for
the Boussinesq equations though it is likely that 
these results can also be obtained through 
other methods.
In the latter case, it is also shown that the associated power series 
in the Borel plane has a radius of convergence independent of size 
of initial data and forcing when initial data and 
forcing have a fixed number of Fourier modes; 
this is useful in computing the solution in the Borel plane.  

\section{Main Results}\label{sec1}

We first write the equations as integral equations in time in
Fourier space. We denote the Fourier transform operator by $\mathcal{F}$, the Fourier transform of $f$ by $\hat{f}$, and $\hat{*}$ the Fourier convolution. As usual, a repeated index $j$ denotes the sum over $j$ from $1$
to $d$. 
$P_k$ is the Fourier transform of the Hodge projection and has the representation
\begin{equation*}
P_k\equiv \left( 1-\frac{k(k\cdot)}{|k|^2}\right).
\end{equation*}
Formal derivation\footnote{While at this stage derivation is formal, 
in the space of functions where existence
is proved, it will become clear that the integral and differential formulations are equivalent.} based on inversion
of the heat operator in Fourier space in (\ref{B}) leads to the following integral equations:
\begin{align}\label{FB2}
\hat{u}(k,t)=- \int_0^t& e^{-\nu|k|^2(t-\tau)}\left(i k_j P_k[\hat u_j\hat *\hat u-a e_2\hat{\Theta}](k,\tau)-\hat f(k)\right)d\tau\\ \nonumber
&+e^{-\nu|k|^2t}\hat u_0(k)\\ \nonumber
\hat{\Theta}(k,t)=- \int_0^t& e^{-\mu|k|^2(t-\tau)}\left( i k_j[\hat u_j\hat *\hat{\Theta}](k,\tau)\right)d\tau+e^{-\mu|k|^2t} \hat{\Theta}_0(k).
\end{align}

\begin{Definition}\label{gamma,beta} We introduce the norm $||\cdot||_{\gamma,\beta}$ for some $\beta \geq 0$ 
and $\gamma >d$:
\begin{equation*}
|| \hat f ||_{\gamma,\beta}=\sup_{k\in \mathbb{R}^d} (1+|k|)^{\gamma}e^{\beta |k|}|\hat f(k)|, \textnormal{ where }\hat{f}(k)=\mathcal{F}[f(\cdot)](k).
\end{equation*}
\end{Definition}
\begin{Definition}\label{L1}
We also use the space $L^1\cap L^{\infty}$ with the norm defined by 
\begin{equation*}
||\hat{f}||_{L^1\cap L^{\infty}}=\max\left\{\int_{\mathbb{R}^d} |\hat{f}(k)|dk, \sup_{k\in\mathbb{R}^d}|\hat{f}(k)|\right\}.
\end{equation*}
\end{Definition}

In cases when results hold either 
for $\| \cdot \|_{\gamma, \beta}$ or  
$\| \cdot \|_{L^1\cap L^{\infty}}$ norm,  
we will use $||\cdot||_{N}$ for brevity of notation. 

We assume 
$||(1+|k|)^2( \hat{u}_0,\hat{\Theta}_0)||_{N}<\infty$
and $||\hat{f}||_{N}<\infty$ in what follows. 
If $\| \cdot \|_N= \| \cdot \|_{\gamma,\beta}$ and $\beta>0$ then the initial condition and forcing are 
real analytic in $x$ in a strip of width at least $\beta$.
 
\begin{Theorem}\label{existence} (Boussinesq Existence and Uniqueness) 

If $\|(1+|\cdot|)^2(\hat{u}_0, \hat{\Theta}_0)\|_{N}
<\infty$ and $||\hat{f}||_{N}<\infty$, then the following statements hold:

i) The Boussinesq equation (\ref{FB2}) has a solution 
$(\hat{u}, \hat{\Theta})(k,t)$ such that $\|(\hat{u}, \hat{\Theta})(\cdot,t)\|_{N}<\infty$ 
for $\Re \frac{1}{t} >\omega$ for $\omega$ sufficiently large\footnote{ 
$\omega$ is large enough so that (\ref{1-2.39}) in the ensuing holds, 
where $(\hat{u}_1,\hat{\Theta}_1)$, defined in (\ref{u1}), depends on the initial data and forcing}. 

ii) The solution has the Laplace transform representation
\begin{equation}\label{Laplacetransform1}
(\hat{u}, \hat{\Theta})(k,t)=(\hat{u}_0, \hat{\Theta}_0)(k)+\int_0^{\infty}(\hat{H}, \hat{S})(k,p)e^{-p/t}dp
\end{equation}
where $(\hat{H},\hat{S})$ is the unique solution to a set of integral equations 
in the space where 
$\|(\hat{H},\hat{S})(\cdot,p)\|_{N}   e^{-\omega p} \in L^1(0,\infty)$. 
The corresponding $(u, \Theta)(x,t)=\mathcal{F}^{-1}[(\hat{u},\hat{\Theta})(k)](x,t)$ is 
analytic in $t$ for $\Re \frac{1}{t}>\omega$ and for $\beta>0$ 
is analytic in $x$ for any $t \in \left [0, \frac{1}{\omega} \right ) $ 
in a strip of width $\beta$, where
initial data and forcing are analytic.

iii) Further, this solution satisfies
$\|(1+|\cdot|)^2(\hat{u}, \hat{\Theta})(\cdot,t)\|_{N} < \infty$ for $t\in (0,\omega^{-1})$,
implying at least the same regularity as initial conditions.
Moreover, $(u, \Theta)(x,t)$ solves (\ref{B}) 
and is the unique Boussinesq solution in $L^{\infty}(0,T,H^2(\mathbb{R}^d))$ when
initial data and forcing in Fourier space satisfy given assumptions. 

iv) A sufficient condition for 
global existence of smooth solution is that
$e^{-\omega p}\|(\hat{H},\hat{S})(\cdot,p)\|_{N}\in L^1(0,\infty)$ for any $\omega> 0$.
\end{Theorem}

\begin{Remark}
{\rm 
If instead we assume $\|(1+|\cdot|^2)(\hat{u}_0,\hat{\Theta}_0)\|_{L^1(\mathbb{R}^d)}<\infty$ and $\|\hat f\|_{L^1(\mathbb{R}^d)}<\infty$, then we have a unique solution to (\ref{FB2}) 
for which $||(\hat{u},\hat{\Theta})||_{L^1(\mathbb{R}^d)}<\infty$ for $t\in (0,\omega^{-1})$. 
Using the arguments of Lemma \ref{2.11.} for $L^1 (\mathbb{R}^d)$ norm alone, the solution is shown to be in the space where $\| (1+ |\cdot|^2) (\hat{u} (\cdot, t), \hat{\Theta} (\cdot, t) \|_{L^1 (\mathbb{R}^d)}$ is finite for $t \in [0,T]$ and solves (\ref{B}) as well. What is not known is 
whether the corresponding $(u, \Theta)$ in the physical $x$-space 
is in $L^\infty(0,T,H^2(\mathbb{R}^d))$.
}
\end{Remark}

\begin{Remark}{\rm The guaranteed existence time $T=\omega^{-1}$ depends 
on $\|(1+|\cdot|)^2(\hat{u}_0,\hat{\Theta}_0)(\cdot)\|_{N}$ and $\| {\hat f} \|_N$. This condition 
is likely to be weakened using an accelerated version of the Borel transform as in \cite{longtime}, {\it i.e.} using
an alternate representation for $n > 1$:
\begin{equation}
\label{accelerated}
(\hat{u},\hat{\Theta})(k,t)=(\hat{u}_0,\hat{\Theta}_0)(k)+\int_0^{\infty}(\hat{H},\hat{S})(k,q)e^{-q/(t^n)}dq
\end{equation}
Further, we expect to prove, that in the periodic case ($x \in \mathbb{T}^d$) without forcing, for any specific
initial condition, global solutions of the PDE implies that 
there exists $n$ sufficiently large so that
$\omega$ for the 
associated integral equation solution is arbitrarily small, 
a result already known \cite{longtime} for the 3-d Navier-Stokes.}
\end{Remark}
\begin{Theorem} \label{Borel summability} (Borel Summability) 

i) For $\beta > 0$, {\it i.e.} for 
analytic initial data and forcing, the Boussinesq solution $(u, \Theta)$ is Borel summable in $t^{-1}$, {\it i.e.}
there exists $(H, S)(x,p)$ analytic in $p$ 
in a neighborhood of $
\{ 0 \} \cup \mathbb{R}^{+}$, 
exponentially bounded for large $p$ and analytic in $x$ for $|$Im$~x_j |<\beta$ for 
$ j = 1, \cdots, d$ such that
\begin{equation}\label{Laplacetransform3}
(u, \Theta)(x,t)=(u_0, \Theta_0)(x)+\int_0^{\infty}(H, S)(x,p)e^{-p/t}dp.
\end{equation}
\noindent In particular, as $t\rightarrow 0^+$,
\begin{equation*} 
(u, \Theta)(x,t)\sim (u_0, \Theta_0)(x)+
\sum_{m=1}^{\infty}(u_m, \Theta_m) (x)t^m,
\end{equation*}
\noindent where $|(u_m, \Theta_m)(x)|\leq m! A_0 D_0^m$ with constants
$A_0$ and $D_0$ generally dependent on the initial condition and forcing through Lemma \ref{3.2.}.

ii) Further, if initial data and forcing have a finite
finite number of Fourier modes, then
the solution $(\hat{H},\hat{S})(k,p)$ has a power series in $p$ with
radius of convergence $D_0^{-1}$ is  
independent of the size of initial data and forcing.
\end{Theorem}

\begin{Remark}
{\rm 
In the case $\beta >0$, we do not need the restriction $\gamma>d$. If $||\hat{u}||_{\gamma,\beta}<\infty$, then for $\beta'\in (0,\beta)$ we have for any $n\in \mathbb{N}$,  $||\hat{u}||_{\gamma+d,\beta'}<\infty$.
}
\end{Remark}

\begin{Remark}
{\rm 
When the the radius of convergence $D_0^{-1}$ 
is independent of size of initial condition and forcing, as is
definitely the case for initial conditions and forcing with
finite Fourier modes, the  
solution can be found conveniently on $[0, p_0]$ through a power series. 
More generally, for specific initial conditions and forcing, the
solution in $[0, p_0]$ may be obtained numerically
with rigorous error bounds similar to NSE \cite{longtime}.
In the following Theorem \ref{improved existence}, we obtain
revised estimates on $\omega$ and therefore existence time of
PDE solution, based on integral equation solution on $[0, p_0]$.
}
\end{Remark}

Let $(\hat{H},\hat{S})(k,p)$ be the solution to (\ref{1-2.18}) provided by Lemma \ref{2.8.}. Define
\begin{equation}\label{beginbehave}
(\hat{H}, \hat{S})^{(a)}(k,p)=\begin{cases}(\hat{H}, \hat{S})(k,p) \mbox{ for } p\in (0,p_0]\subset \mathbb{R}^{+}\\
0 \mbox{ otherwise}\end{cases}
\end{equation}
\noindent and
\begin{align}\nonumber\label{smallbehave}
\hat{H}^{(s)}(k,p)&=\frac{i k_j\pi}{2|k|\sqrt{\nu p}}\int_0^{\min(p,2p_0)} \mathcal{G}(z,z')\hat{G}_j^{[1],(a)}(k,p')dp'+2\hat u_1(k)\frac{J_1(2|k|\sqrt{\nu p})}{2|k|\sqrt{\nu p}}\\ 
&\quad+\frac{a\pi}{2|k|\sqrt{\nu p}}\int_0^{\min(p,p_0)} \mathcal{G}(z,z')P_k[e_2\hat{S}^{(a)}(k,p')]dp'\\ \nonumber
\hat{S}^{(s)}(k,p)&=\frac{i k_j\pi}{2|k|\sqrt{\mu p}}\int_0^{\min(p,2p_0)} \mathcal{G}(\zeta ,\zeta ')\hat{G}_j^{[2],(a)}(k,p')dp'+2\hat \Theta_1(k)\frac{J_1(2|k|\sqrt{\mu p})}{2|k|\sqrt{\mu p}}
\end{align}
\noindent where 
\begin{align*}
\hat{G}_j^{[1],(a)}(k,p)&=-P_k[\hat{u}_{0,j}\hat{*}\hat{H}^{(a)}+\hat{H}_j^{(a)}\hat{*}\hat{u}_0+\hat{H}_j^{(a)}\, ^{\ast}_{\ast}\hat{S}^{(a)}]\\ 
\hat{G}_j^{[2],(a)}(k,p)&=-[\hat{u}_{0,j}\hat{*}\hat{S}^{(a)}+\hat{H}_j^{(a)}\hat{*}\hat{\Theta}_0+\hat{S}_j^{(a)}\, ^{\ast}_{\ast}\hat{S}^{(a)}].
\end{align*}
Notice if $(\hat{H}, \hat{S})^{(a)}(k,p)$ is known, then $\hat{H}^{(s)}(k,p)$, $\hat{S}^{(s)}(k,p)$, $G_j^{[1],(a)}(k,p)$, and $G_j^{[2],(a)}(k,p)$ are also known functions. Also, recall $\hat{u}_1$ and $\hat{\Theta}_1$ are quantities based on the initial condition and forcing given in (\ref{u1}).

\begin{Theorem}{(Revised Exponential Estimates)}\label{improved existence}. For
some $\omega_0 \ge 0$, assume $\epsilon_1$, $B_3$ and $b$ are
functionals of the forcing $f$, initial condition $(\hat{u}_0,\hat{\Theta}_0)$, and 
the solution $({\hat H}, {\hat S} )$ to the
set of integral equations (\ref{N}) on
a finite interval $[0, p_0]$,
determined from the relations:
\begin{equation}\label{13.10}
b=\omega_0 \int_{p_0}^{\infty}e^{-\omega_0 p}||(\hat{H}, \hat{S})^{(s)}(\cdot,p)||_{N}dp
\end{equation}
\begin{equation}\label{13.11}
\epsilon_1=\mathcal{B}_1+\mathcal{B}_4+\int_0^{p_0}e^{-\omega_0 p}\mathcal{B}_2(p)dp,
\end{equation}
\noindent where
\begin{multline*}
\mathcal{B}_0(k)=C_0\sup_{p_0\leq p'\leq p}|\mathcal{G}(z,z')/z|,\hspace{.5 in} \mathcal{B}_1=2\sup_{k\in\mathbb{R}^d}|k|\mathcal{B}_0(k)||(\hat{u}_0, \hat{\Theta}_0)||_{N},\\ 
\mathcal{B}_2=2\sup_{k\in\mathbb{R}^d}|k|\mathcal{B}_0(k)||(\hat{H}, \hat{S})^{(a)}(\cdot,p)||_{N},\,\, \mathcal{B}_3=\sup_{k\in\mathbb{R}^d}|k|\mathcal{B}_0(k), \,\,\mathcal{B}_4=a\sup_{k\in\mathbb{R}^d}\mathcal{B}_0(k).
\end{multline*}
Then, over an extended interval $\mathbb{R}^+$, the solution
satisfies the relation
\begin{equation*}
\left \| \left ( {\hat H} (\cdot, p), {\hat S} (\cdot, p) \right ) 
\right \|_{N}e^{-\omega p} \in L^1 \left (0, \infty \right )
\end{equation*}
for any $\omega\geq \omega_0$ satisfying 
\begin{equation*} 
\omega>\epsilon_1+2\sqrt{\mathcal{B}_3 b}. 
\end{equation*}

\end{Theorem}
\begin{Remark}{\rm The implication of the above theorem is that 
if solution $(\hat{H},\hat{S})$, 
restricted to $[0,p_0]$ is known, through computation of power series in $p$ or otherwise, 
and if the corresponding functionals $\epsilon$ and $\mathcal{B}_3 b$ are
small, as is the case for sufficiently rapidly decaying $({\hat H}, {\hat S})$
over a large enough interval $[0, p_0]$, then     
existence for Boussinesq PDE solution in a long interval $(0,\omega^{-1})$ 
is guaranteed. 
It is to be noted that rigorous error control of computed solution in $[0, p_0)$ is
expected as for 3-d NSE \cite{longtime}; this  
leads to a revised bound on $\omega$ that can translate to
a longer existence time.
}
\end{Remark}

\section{Local Existence and Uniqueness of Solution}

\subsection{Formulation of Integral Equation: Borel Transform}\label{sec2}

Our goal is to take the Borel transform and create equivalent 
integral equations. To ensure smallness in $t$ for small $t$ and 
avoid dealing with delta distribution in Borel transform, 
it is convenient to define $\hat{h}$ and $\hat{w}$ so that 
\begin{equation*}
(\hat{u},\hat{\Theta})(k,t)=(\hat{u}_0,\hat{\Theta}_0)(k)+(\hat{h},\hat{s})(k,t).
\end{equation*}

\noindent For (\ref{FB2}), we define 
\begin{equation}\label{1-2.7}
\hat{g}^{[1]}_j:= P_k[\hat{h}_j\hat{*}\hat{h}+\hat{h}_j\hat{*}\hat{u}_0+
\hat{u}_{0,j}\hat{*}\hat{h}]\textnormal{ and }\hat{g}^{[2]}_j:=[\hat{h}_j\hat{*}\hat{s}+\hat{h}_j\hat{*}\hat{\Theta}_0+
\hat{u}_{0,j}\hat{*}\hat{s}]
\end{equation}
and
\begin{align}\label{u1}
\hat{u}_1(k)&:=-\nu|k|^2\hat{u}_0-i k_j P_k[\hat{u}_{0,j}\hat{*}\hat{u}_0]+a P_k[e_2\hat{\Theta}_0]+\hat{f}\\ \nonumber
\hat{\Theta}_1(k)&:=-\mu|k|^2\hat{\Theta}_0-i k_j(\hat{u}_{0,j}\hat{*}\hat{\Theta}_0).
\end{align}

Using these in (\ref{FB2}),
we obtain integral equations:
\begin{align}\label{IE-t-1}
\hat{h}(k,t)&=\int_0^te^{-\nu|k|^2(t-s')}\left(-ik_j\hat{g}^{[1]}_j-P_k[ae_2\hat{s}]\right)(k,s')ds'+\left(\frac{1-e^{-\nu|k|^2t}}{\nu|k|^2}\right)\hat{u}_1\\ \nonumber
\hat{s}(k,t)&=-ik_j\int_0^te^{-\mu|k|^2(t-s')}\hat{g}^{[2]}_j(k,s')ds'+\left(\frac{1-e^{-\mu|k|^2t}}{\mu|k|^2}\right)\hat{\Theta}_1.
\end{align}
We seek a solution as a Laplace transform, 
\begin{equation*} 
(\hat{h},\hat{s})(k,t)=\int_0^{\infty}\left(\hat{H}, \hat{S}\right)(k,p)e^{-p/t}dp.
\end{equation*}
\noindent With this goal, 
we take the formal\footnote[1]{While the derivation of the integral equation in $p$ is formal, 
we prove later (Lemma \ref{2.9.}) that the unique solution to the 
integral equation in the Borel plane generates a solution to the Boussinesq equation through Laplace transform.}  
inverse Laplace transform in $1/t$:
\begin{equation*}
[\mathcal{L}^{-1}f](p)=\frac{1}{2\pi i}\int_{c-i\infty}^{c+i\infty}f(s)e^{sp}ds,
\end{equation*}
where $c$ is chosen so that for $\mathrm{Re} \,s\geq c$, $f$ is analytic and has suitable asymptotic decay. We define 
\begin{equation}\label{Chap2H}
\mathcal{H}^{(\nu)}(p,p',k):=\int_{p'/p}^1\left\{\frac{1}{2\pi i}\int_{c-i\infty}^{c+i\infty}\tau^{-1}exp[-\nu|k|^2\tau^{-1}(1-s)+(p-p's^{-1})\tau]d\tau\right\}ds.
\end{equation}
Then (\ref{IE-t-1}) becomes
\begin{align}\label{IE-p-1}
\hat{H}(k,p)&=\int_0^p\mathcal{H}^{(\nu)}(p,p',k)\left(-ik_j\hat{G}^{[1]}_j(k,p')dp'+P_k[ae_2\hat{S}](k,p')\right)dp'\\ \nonumber
&\qquad\qquad +\hat{u}_1(k)\mathcal{L}^{-1}\left(\frac{1-e^{-\nu|k|^2t}}{\nu|k|^2}\right)(p)\\ \nonumber
\hat{S}(k,p)&=-ik_j\int_0^p\mathcal{H}^{(\mu)}(p,p',k)\hat{G}^{[2]}_j(k,p')dp'+\hat{\Theta}_1(k)\mathcal{L}^{-1}\left(\frac{1-e^{-\mu|k|^2t}}{\mu|k|^2}\right)(p).
\end{align}
\noindent In the above, $\hat{G}_j^{1,2}=\mathcal{L}^{-1}[g_j^{1,2}]$. Specifically,
\begin{equation}\label{G_j}
 \hat{G}^{[1]}_j=P_k[\hat u_{0,j}\hat{*}\hat H+\hat H_j\hat{*}\hat u_0+\hat H_j \, ^*_* \hat H]\,\textnormal {and }\hat{G}^{[2]}_j=[\hat u_{0,j}\hat{*}\hat S+\hat H_j\hat{*}\hat{\Theta}_0+\hat H_j\, ^*_* \hat S]
\end{equation}
where $\, ^* _*$ denotes the Laplace convolution followed by Fourier convolution (order is unimportant). 
We now make the observation that our kernel $\mathcal{H}^{(\nu)}(p,p',k)$ has a representation in terms of Bessel functions. Namely,
\begin{equation*}
\mathcal{H}^{(\nu)}(p,p',k)=\frac{\pi}{z}\mathcal{G}(z,z'):=\frac{\pi z'}{z}\left\{-J_1(z)Y_1(z')+Y_1(z)J_1(z')\right\}
\end{equation*}
where $J_1$ and $Y_1$ are the Bessel functions of order 1, $z=2|k|\sqrt{\nu p}$, and $z'=2|k|\sqrt{\nu p'}$. In similar spirit, we have
\begin{equation*}
\frac{2J_1(z)}{z}=\mathcal{L}^{-1}\left(\frac{1-e^{-\nu|k|^2\tau^{-1}}}{\nu|k|^2}\right)(p).
\end{equation*}
These assertions are proved in the appendix in \ref{kernel} and \ref{U_0}.
Thus, our integral Boussinesq equation becomes 
\begin{align}\label{1-2.18}
\hat{H}(k,p)=&\pi \int_0^p \frac{\mathcal{G}(z,z')}{z}\left(i k_j\hat{G}_j^{[1]}(k,p')+aP_k[e_2\hat S(k,p')]\right)dp'+2\hat u_1(k)\frac{J_1(z)}{z}\\ \nonumber
\hat{S}(k,p)=&\frac{i k_j\pi}{2|k|\sqrt{\mu p}}\int_0^p \mathcal{G}(\zeta ,\zeta ')\hat{G}_j^{[2]}(k,p')dp'+2\hat \Theta_1(k)\frac{J_1(\zeta)}{\zeta},
\end{align}
where $\zeta=2|k|\sqrt{\mu p}$, and $\zeta'=2|k|\sqrt{\mu p'}$. Abstractly, we may write the set of equations (\ref{1-2.18}) as
\begin{equation}\label{N}
(\hat{H},\hat{S})(k,p)=\mathcal{N}[(\hat{H},\hat{S})](k,p).
\end{equation}

\begin{Remark}{\rm 
By properties of Bessel functions $|\mathcal{G}(z,z')|$ is bounded for all real nonnegative $z'\leq z$. (The approximate bound is $0.6$, see \cite{smalltime}). The asymptotic properties of Bessel functions for small $z$ also show $|\mathcal{G}(z,z')/z|$ is bounded for all real nonnegative $z'\leq z$.
}
\end{Remark}

To prove Theorem \ref{existence}, we will show $\mathcal{N}$ is contractive in a suitable space, 
so $(\hat H,\hat S)$ is Laplace transformable in $1/t$. 
Then from Lemma \ref{2.9.}
\begin{equation*}
(\hat u,\hat{\Theta})(k,t)=(\hat u_0,\hat{\Theta}_0)(k)+\int_0^{\infty}(\hat H,\hat{S})(k,p)e^{-p/t}dp
\end{equation*}
satisfies ($\ref{FB2}$) for $\Re \left ( 1/t \right ) $ large enough.
Furthermore, we show $(u, \Theta)(x,t)=\mathcal{F}^{-1}[(\hat{u},\hat{\Theta})(\cdot,t)](x)$ 
is a classical solution to the Boussinesq problem. 

\subsection{Norms in p}\label{sec3}
Recall the norm $|| \cdot ||_{N}$ in $k$ is either the $(\gamma,\beta)$ norm given in Definition \ref{gamma,beta} for some $\beta \geq 0$ and $\gamma >d$ or the $L^1\cap L^{\infty}$ norm.
\begin{Definition}
For $\alpha \geq 1$, we define 
\begin{equation*}
||\hat f||^{(\alpha)}=\sup_{p\geq 0}(1+p^2)e^{-\alpha p}||\hat{f}(\cdot, p)||_{N}.
\end{equation*}
\end{Definition}

\begin{Definition}
We define $\mathcal{A}^{\alpha}$ to be the Banach space of continuous function of $(k,p)$ for $k\in \mathbb{R}^d$ and $p\in \mathbb{R}^+$ for which $||\cdot ||^{\alpha}$ is finite. In similar spirit, we define the space $\mathcal{A}_1^{\alpha}$ of 
locally integrable functions for $p\in [0,L)$, and continuous in $k$ such that
\begin{equation*}
||\hat f||_1^{\alpha}=\int_0^Le^{-\alpha p}||\hat f(\cdot, p)||_{N}dp < \infty.
\end{equation*}
\end{Definition}
\begin{Definition}
Finally, we also define $\mathcal{A}_L^{\alpha}$ to be the Banach space of continuous functions in $(k,p)$ for $k$ in $\mathbb{R}^d$ and $p\in [0,L]$ such that
\begin{equation*}
||\hat f||^{\infty}_L=\sup_{p\in[0,L]}||\hat f(\cdot,p)||_{N}<\infty.
\end{equation*}
\end{Definition}

\subsection{Existence of a Solution in Dual Variable}\label{sec4}
We need some preliminary lemmas. Recall, $d=2$ or $d=3$ denotes the dimension in $x$ or its dual $k$. Often constants appearing in subalgebra bounds will depend on dimension. We will explicitly state the dependence when defining them and suppress the dependence elsewhere.

\begin{Lemma}
If $||\hat{v}||_{\gamma, \beta}$ and $||\hat{w}||_{\gamma, \beta}<\infty$ for $\gamma >d$ and $k\in \mathbb{R}^d$, then
\begin{equation*}
||\hat{v}\hat{*}\hat{w}||_{\gamma,\beta}\leq \tilde{C}_0(d)||\hat{v}||_{\gamma, \beta}||\hat{w}||_{\gamma, \beta},
\end{equation*}
\noindent where 
\begin{align*}
\tilde{C}_0(2)=2^{\gamma+1}\int_{k'\in \mathbb{R}^2}\frac{1}{(1+|k'|)^{\gamma}}dk'=\frac{\pi 2^{\gamma +2}}{(\gamma -1)(\gamma -2)}  \textnormal{ and }\\
\tilde{C}_0(3)=2^{\gamma+1}\int_{k'\in \mathbb{R}^3}\frac{1}{(1+|k'|)^{\gamma}}dk'=\frac{\pi 2^{\gamma +4}}{(\gamma -1)(\gamma -2)(\gamma -3)}.
\end{align*} 
\end{Lemma}

\begin{pf} The $d=3$ case
can be found in \cite{smalltime} and 
the $d=2$ case is basically the same. From the definition of $||\cdot||_{\gamma,\beta}$ and the fact that $e^{-\beta(|k'|+|k-k'|)}\leq e^{-\beta|k|}$, we have 
\begin{equation*}
|\hat{v}\hat{*}\hat{w}|\leq e^{-\beta|k|}||\hat{v}||_{\gamma,\beta}||\hat{w}||_{\gamma,\beta}\int_{k'\in\mathbb{R}^2}(1+|k'|)^{-\gamma}(1+|k-k'|)^{-\gamma}dk'.
\end{equation*}
\noindent Split the integral into two domains $|k'|\leq|k|/2$ and its complement to show
\begin{align*}
\int_{k'\in\mathbb{R}^2}\frac{1}{(1+|k'|)^{\gamma}(1+|k-k'|)^{\gamma}}dk'&\leq \frac{2^{\gamma+1}}{(1+|k|)^{\gamma}}\int_{k'\in\mathbb{R}^2}\frac{1}{(1+|k'|)^{\gamma}}dk'\\ 
&=\frac{2^{\gamma+2}\pi}{(1+|k|)^{\gamma}(\gamma-1)(\gamma-2)},
\end{align*} 
where polar coordinates and integration by parts are used to evaluate the last integral.
\end{pf}

\begin{Corollary}\label{2.2.} If $||\hat{v}||_{N}$, $||\hat{w}||_{N}<\infty$, then for $C_0=C_0(d)$ chosen such that $C_0=\tilde{C}_0 $ for $N = (\gamma, \beta) $, $\gamma>d$ and $C_0=1$ for $N=L^1\cap L^{\infty}$, we have
\begin{equation*}
||\hat{v}\hat{*}\hat{w}||_N\leq C_0||\hat{v}||_N||\hat{w}||_N.
\end{equation*}
\end{Corollary}

\begin{Lemma}\label{2.3.} Also, notice that
\begin{equation*}
\left\|\left(P_k(\hat{f}), P_k(\hat{g})\right)\right\|_{N}\leq ||(\hat{f} , \hat{g})||_{N}
\end{equation*}
\end{Lemma}
\begin{pf} $P_k$ is the projection of a vector onto $k^{\bot}$.
\end{pf}

\begin{Lemma}\label{2.4.} With $C_0$ as defined in Corollary \ref{2.2.}, appropriately modified for $d=2$ or $3$, and constants
\begin{equation*}
C_2=\frac{\pi C_0}{\min (\sqrt{\nu}, \sqrt{\mu})} \sup_{z\in\mathbb{R}^+,0\leq z'\leq z}|\mathcal{G}(z,z')|\textnormal{ and }C_3=\pi a\sup_{z\in\mathbb{R}^+,0\leq z'\leq z}\left|\frac{\mathcal{G}(z,z')}{z}\right|,
\end{equation*} 
we have the following bounds on the norm in $k$ for the operator $\mathcal{N}$ defined in (\ref{N}). Let $\phi:=(\hat H,\hat S)$. Then
\begin{multline}\label{1-2.23}
||\mathcal{N}[\phi(\cdot, p)]||_{N} \leq  \frac{C_2}{\sqrt{p}}\int_0^p\left( ||\phi (\cdot, p')||_{N}*||\phi(\cdot,p')||_{N}\right. \hspace{2 in}\\ 
+\left. ||(\hat{u}_0, \hat \Theta_0)||_{N}||\phi(\cdot,p')||_{N} \right) dp'
+||(\hat u_1, \hat \Theta_1)||_{N}+C_3\int_0^p||\hat S(\cdot,p')||_{N}dp'
\end{multline}
and
\begin{align}\label{1-2.24}
||\mathcal{N}&[\phi^{[1]}](\cdot, p)-\mathcal{N}[\phi^{[2]}](\cdot, p)||_{N} \leq \frac{C_2}{\sqrt{p}}\int_0^p \left( ||\phi^{[1]} (\cdot, p')||_{N}+||\phi^{[2]}(\cdot, p')||_{N}\right)\\ \nonumber
& *\left\|\phi^{[1]}-\phi^{[2]}(\cdot,p')\right\|_{N} +||(\hat u_0,\hat \Theta_0)||_{N}||\phi^{[1]}-\phi^{[2]}(\cdot,p')||_{N} dp'\\ \nonumber
&+C_3\int_0^p||\hat S^{[1]}-\hat S^{[2]}(\cdot, p')||_{N}dp'
\end{align}
\end{Lemma}

\begin{pf} From \cite{handbook}, $|J_1(z)/z|\leq 1/2$ for $z\in \mathbb{R}^+$ and
\begin{equation*}
\left\|2\left(\hat u_1(k)\frac{J_1(z)}{z}, \hat \Theta_1(k)\frac{J_1(\zeta)}{\zeta}\right)\right\|_{N}\leq ||(\hat u_1, \hat {\Theta}_1)||_{N}.
\end{equation*}

\noindent From Corollary \ref{2.2.}, we have
\begin{multline*}
|||\hat u_{0}|\hat{*}(\hat H, \hat{S})+|\hat H|\hat{*}(\hat u_0, \hat{\Theta}_0)+|\hat H| \, ^*_* (\hat H,\hat{S})||_{N}\leq \hspace{2 in}\\ 
\left[2C_0||(\hat{u}_0,\hat{\Theta}_0)||_{N}||(\hat{H},\hat{S})(\cdot, p)||_{N}+C_0||\hat{H}(\cdot, p)||_{N}*||(\hat{H},\hat{S})(\cdot, p)||_{N}\right].
\end{multline*}
\noindent Then using Lemma \ref{2.3.} and Schwartz inequality, we obtain
\begin{equation*}
||k_j(\hat{G}_j^{[1]}, \hat{G}_j^{[2]})||_{N}\leq 2C_0|k|\left(||\phi(\cdot, p')||_{N}*||\phi(\cdot,p')||_{N}+ ||(\hat{u}_0, \hat {\Theta}_0)||_{N}||\phi(\cdot,p')||_{N} \right).
\end{equation*}
\noindent Now (\ref{1-2.23}) follows.  To obtain (\ref{1-2.24}) notice that
\begin{equation}\label{Chap2num}
\hat H_j^{[1]} \, ^*_*\phi^{[1]}-\hat H_j^{[2]} \, ^*_*\phi^{[2]}=\hat H_j^{[1]} \, ^*_* \left(\phi^{[1]}-\phi^{[2]}\right)+(\hat H_j^{[1]}-\hat H_j^{[2]}) \, ^*_*\phi^{[2]}.
\end{equation}
\noindent From (\ref{Chap2num}) we get 
\begin{equation*} 
\left\|\hat H_j^{[1]} \, ^*_*\phi^{[1]}-\hat H_j^{[2]} \, ^*_*\phi^{[2]}\right\|_{N}\leq C_0 \left\|\phi^{[1]}-\phi^{[2]}\right\|_{N} *\left( ||\phi^{[1]}||_{N} +||\phi^{[2]}||_{N} \right).
\end{equation*}
Combining this bound and using Lemma \ref{2.3.} as in the first part of the proof, we get (\ref{1-2.24}). 
\end{pf}
\begin{Lemma}\label{2.6.} For $\hat{f},\hat{g} \in \mathcal{A}^{\alpha}, \mathcal{A}^{\alpha}_1$ or $\mathcal{A}^{\infty}_L$
\begin{align*}
||\hat{f}\,^*_*\hat{g}||^{(\alpha)}&\leq M_0 C_0 ||\hat{f}||^{(\alpha)}||\hat{g}||^{(\alpha)}\\ 
||\hat{f}\,^*_*\hat{g}||^{\alpha}_1&\leq C_0||\hat{f}||^{\alpha}_1||\hat{g}||^{\alpha}_1\\ 
||\hat{f}\,^*_*\hat{g}||^{\infty}_L&\leq L C_0 ||\hat{f}||^{\infty}_L||\hat{g}||^{\infty}_L,
\end{align*} 
\noindent where $M_0\approx 3.76\cdots$ is large enough so
\begin{equation*}
\int_0^p\frac{(1+p^2)ds}{(1+s^2)(1+(p-s)^2)}\leq M_0.
\end{equation*}
\end{Lemma}
\noindent This means the Banach spaces listed in the norms section form subalgebras under the operation $\,^*_*$. The properties listed are independent of dimension except for a change in $C_0$ showing up due to the Fourier convolution. The proof is in \cite{smalltime}. The basic idea is that $k$ and $p$ act separately in the norm. So, 
we need only consider how the $p$ portion of the norm effects $\int_0^p u(p)v(p-s)ds$.

The following lemma expands the bounds in Lemma \ref{2.4.} to bounds in $p$ in some of our other norms. 

\begin{Lemma}\label{2.7.} 
Let $\phi:=(\hat H,\hat S)$. On $\mathcal{A}_1^{\alpha}$, the operator $\mathcal{N}$ satisfy the following inequalities 
\begin{multline}\label{1-2.30}
||\mathcal{N}[\phi]||_1^{\alpha}\leq C_2 \sqrt{\pi}\alpha^{-1/2}\left\{ (||\phi||_1^{\alpha})^2+||(\hat u_0, \hat{\Theta}_0)||_{N}||\phi||_1^{\alpha}\right\} \\ 
+\alpha ^{-1}||(\hat u_1, \hat \Theta_1)||_{N}+\alpha^{-1}C_3||\hat S||_1^{\alpha}
\end{multline}
and
\begin{multline}\label{1-2.31}
||\mathcal{N}(\phi^{[1]})-\mathcal{N}(\phi^{[2]})||_1^{\alpha}\leq C_2 \sqrt{\pi}\alpha^{-1/2}\left\{ \left(||\phi^{[1]}||_1^{\alpha}+||\phi^{[2]}||_1^{\alpha}\right)\right.\\
\left(||\phi^{[1]}-\phi^{[2]}||_1^{\alpha}\right) +\left. ||(\hat u_0, \hat{\Theta}_0)||_{N}||\phi^{[1]}-\phi^{[2]}||_1^{\alpha} \right\} +\alpha^{-1}C_3||\hat S^{[1]}-\hat S^{[2]}||_1^{\alpha},
\end{multline}
\noindent Similarly, for $\mathcal{A}_L^{\infty}$, we have
\begin{multline}\label{1-2.32}
||\mathcal{N}[\phi]||_L^{\infty}\leq C_2 \sqrt{L}\left\{ L(||\phi||_L^{\infty})^2+||(\hat u_0, \hat{\Theta}_0)||_{N}||\phi||_L^{\infty}\right\} \hspace{1.05 in}\\
+||(\hat u_1, \hat \Theta_1)||_{N}+LC_3||\hat S||_L^{\infty}
\end{multline}
and
\begin{multline} \label{1-2.33}
||\mathcal{N}[\phi^{[1]}]-\mathcal{N}[\phi^{[2]}]||_L^{\infty}\leq C_2 \sqrt{L}\left\{ L\left(||\phi^{[1]}||_L^{\infty}+||\phi^{[2]}||_L^{\infty}\right)\right.\\
\left(||\phi^{[1]}-\phi^{[2]}||_L^{\infty}\right)+\left. ||(\hat u_0, \hat{\Theta}_0)||_{N}||\phi^{[1]}-\phi^{[2]}||_L^{\infty} \right\} + LC_3||\hat S^{[1]}-\hat S^{[2]}||_L^{\infty},
\end{multline}
\end{Lemma}

\begin{pf} For the space $\mathcal{A}_1^{\alpha}$ and any $L>0$, we note that
\begin{equation*}
\int_0^Le^{-\alpha p}||(\hat u_1,\hat \Theta_1)||_{N}dp\leq \alpha^{-1}||(\hat u_1,\hat \Theta_1)||_{N}
\end{equation*}
\noindent and
\begin{equation*}
\int_0^L e^{-\alpha p}p^{-1/2}dp\leq \Gamma \left( \frac{1}{2} \right) \alpha^{-1/2}=\sqrt{\pi}\alpha^{-1/2}.
\end{equation*}
\noindent We further notice that for $y(p')\geq 0$, we have
\begin{multline}
\int_0^L e^{-\alpha p}p^{-1/2}\left( \int_0^p y(p')dp'\right) dp=\int_0^L y(p')e^{-\alpha p'}\left( \int_{p'}^L e^{-\alpha (p-p')}p^{-1/2}dp\right) dp'\\ \label{2.34}
\leq \int_0^L y(p')e^{-\alpha p'}\left( \int_{0}^L e^{-\alpha s}s^{-1/2}ds\right) dp'\leq \int_0^L y(p')e^{-\alpha p'}\sqrt{\pi}\alpha^{-1/2}dp'.
\end{multline}
\noindent Similarly,
\begin{equation*}
\int_0^L e^{-\alpha p}\left( \int_0^p ||\hat S(\cdot,p')||_{N}dp'\right) dp\leq \alpha ^{-1}||\hat S||_1^{\alpha}.
\end{equation*}
\noindent Then, using (\ref{2.34}) in ($\ref{1-2.23}$) and the idea in Lemma \ref{2.6.} that $\int_0^p e^{-\alpha p}[(||g||_{N}*||h||_{N})(p)]dp\leq ||g||_1^{\alpha}||h||_1^{\alpha}$, we have 
\begin{multline*}
\int_0^Le^{-\alpha p}||\mathcal{N}(\hat H, \hat S)||_{N} dp\leq C_2 \sqrt{\pi}\alpha^{-1/2}\left\{ (||(\hat H, \hat S)||_1^{\alpha})^2\right. \hspace{1 in}\\
+\left. ||(\hat u_0, \hat{\Theta}_0)||_{N}||(\hat H, \hat S)||_1^{\alpha} \right\} +\alpha ^{-1}||(\hat u_1, \hat \Theta_1)||_{N}+\alpha^{-1}C_3||\hat S||_1^{\alpha}.
\end{multline*}

\noindent This proves (\ref{1-2.30}). Further, from ($\ref{1-2.24}$), it also follows that
\begin{multline*}
\int_0^L e^{-\alpha p}||\mathcal{N}(\phi^{[1]})-\mathcal{N}(\phi^{[2]})(\cdot,p)||_{N}dp\leq C_2 \sqrt{\pi}\alpha^{-1/2}\left\{ \left(||\phi^{[1]}||_1^{\alpha}+||\phi^{[2]}||_1^{\alpha}\right)\right.\\
\left\|\phi^{[1]}-\phi^{[2]}\right\|_1^{\alpha} +\left. ||(\hat u_0, \hat{\Theta}_0)||_{N}\left\|\phi^{[1]}-\phi^{[2]}\right\|_1^{\alpha} \right\} +\alpha^{-1}C_3||\hat S^{[1]}-\hat S^{[2]}||_1^{\alpha}.
\end{multline*}

\noindent This proves (\ref{1-2.31}).

Now, we consider $\mathcal{A}_L^{\infty}$. We note that for $p\in[0,L]$, we have
\begin{equation*}
\left| p^{-1/2}\int_0^p y(p')dp'\right| \leq \sup_{p\in [0,L]}|y(p)|\sqrt{L}.
\end{equation*}
\noindent We recall from Lemma \ref{2.6.} that
\begin{equation*}
\left| \int_0^p y_1(s)y_2(p-s)ds\right| \leq L\left( \sup_{p\in [0,L]}|y_1(p)|\right) \left( \sup_{p\in [0,L]}|y_2(p)|\right). 
\end{equation*}
\noindent Taking
\begin{align*}
&y(p)=||\phi(\cdot,p)||_{N}*||\phi(\cdot,p)||_{N}+||(\hat u_0, \hat \Theta_0)||_{N}||\phi(\cdot ,p)||_{N}\\
&\textnormal{and }y_1(p)=y_2(p)=||\phi(\cdot,p)||_{N},
\end{align*}

\noindent ($\ref{1-2.32}$) follows from ($\ref{1-2.23}$). To get the bound in ($\ref{1-2.33}$), we will choose
\begin{multline*}
y(p)= \left(||\phi^{[1]}||_{N}+||\phi^{[2]}||_{N}\right)*\left\|\phi^{[1]}-\phi^{[2]}\right\|_{N} + ||(\hat u_0, \hat{\Theta}_0)||_{N}||\phi^{[1]}-\phi^{[2]}||_{N}, \\
y_1(p)=||\phi^{[1]}||_{N}+||\phi^{[2]}||_{N}, \textnormal{ and }
y_2(p)=\left\|\phi^{[1]}-\phi^{[2]}\right\|_{N} 
\end{multline*}
\noindent now using ($\ref{1-2.24}$) the proof follows. 
\end{pf}

\begin{Lemma}\label{2.8.} Equation ($\ref{1-2.18}$) has a unique solution in $\mathcal{A}_1^{\omega}$ for any $L > 0$ in a ball of size $2\omega^{-1}||(\hat{u}_1,\hat \Theta_1)||_{N}$ for $\omega$ large enough to guarantee
\begin{equation}\label{1-2.39}
2C_2\sqrt{\pi}\omega^{-1/2}\left\{2\omega^{-1} ||(\hat u_1, \hat \Theta_1)||_{N}+||(\hat u_0,\hat \Theta_0)||_{N} +\frac{C_3}{C_2\sqrt{\pi}}\,\omega^{-1/2}\right\}<1
\end{equation}
where $(\hat{u}_1,\hat{\Theta}_1)$ is given in (\ref{u1}). Furthermore, the solution also belongs to $\mathcal{A}_L^{\infty}$ for $L$ small enough to ensure 
\begin{equation}\label{1-2.40}
2C_2L^{1/2}\left\{ 2L||(\hat u_1, \hat \Theta_1)||_{N}+||(\hat u_0,\hat \Theta_0)||_{N}+\frac{C_3}{C_2}L^{1/2}\right\} <1.
\end{equation}
Moreover, $\lim_{p\rightarrow 0^+}(\hat H,\hat S)(k,p)=(\hat u_1, \hat \Theta_1)(k)$. 
\end{Lemma}
\begin{pf} The estimates in Lemma \ref{2.7.} imply that $\mathcal{N}$ maps a ball of radius $2\omega^{-1}||(\hat u_1, \hat{\Theta}_1)||_{N}$ in $\mathcal{A}_1^{\omega}$ into itself and is contractive when $\omega$ is large enough to satisfy ($\ref{1-2.39}$). Similarly, $\mathcal{N}$ maps a ball of size $2||(\hat u_1, \hat{\Theta}_1)||_{N}$ in $\mathcal{A}_L^{\infty}$ into itself and is contractive when $L$ is small enough to satisfy ($\ref{1-2.40}$). Therefore, there is a unique solution to the Boussinesq integral system of equations in the ball. Furthermore, $\mathcal{A}_L^{\infty}\subseteq \mathcal{A}_1^{\alpha}$, so the solutions are in fact one and the same. 

Moreover, applying (\ref{1-2.33}) with $(\hat{H}^{[1]},\hat{S}^{[1]})=(\hat{H},\hat{S})$ and $(\hat{H}^{[2]},\hat{S}^{[2]})=0$, we obtain
\begin{multline*}
\left\|(\hat{H},\hat{S})(k,p)-\left(\hat{u}_1(k)\frac{2J_1(z)}{z},\hat{\Theta}_1(k)\frac{2J_1(\zeta)}{\zeta}\right)\right\|_L^{\infty}\leq\hspace{2 in}\\  C_2L^{1/2}\left\{L(||(\hat{H},\hat{S})||_L^{\infty})^2+||(\hat{u}_0,\hat{\Theta}_0)||_{N}||(\hat{H},\hat{S})||_L^{\infty}\right\}+LC_3||\hat{S}||_L^{\infty}.
\end{multline*}
\noindent Since $||(\hat{H},\hat{S})||_L^{\infty}$ is bounded for small $L$, letting $L\rightarrow 0$,
\begin{equation*}
\left\|(\hat{H},\hat{S})(k,p)-\left(\hat{u}_1(k)\frac{2J_1(z)}{z},\hat{\Theta}_1(k)\frac{2J_1(\zeta)}{\zeta}\right)\right\|_L^{\infty}\rightarrow 0.
\end{equation*}
As $\lim_{z\rightarrow 0}2J_1(z)/z=1$, for fixed $k$,  
$\lim_{p\rightarrow 0}(\hat{H},\hat{S})(k,p)=(\hat{u}_1,\hat{\Theta}_1)(k)$. 
\end{pf}

\subsection{Proof of Local Existence for Boussinesq PDE}\label{sec5}

We have unique solutions to our integral equation, (\ref{IE-p-1}). We show in the following Lemma \ref{2.9.} that the solution's Laplace transform gives a solution to (\ref{FB2}), which is analytic in $t$ for $\Re \frac{1}{t}>\omega$. Lemma \ref{equivalence} below shows that any solution of (\ref{FB2}) with $||(1+|\cdot|)^2(\hat{u},\hat{\Theta})(\cdot, t)||_{N}<\infty$ is inverse Fourier transformable with $(u,\Theta)$ solving (\ref{B}).  Lemma \ref{2.11.} below ensures that $||(1+|\cdot|)^2(\hat{u},\hat{\Theta})(\cdot, t)||_{N}<\infty$. Thus, combining these results, we have $(u,\Theta)(x,t)=\mathcal{F}^{-1}(\hat{u},\hat{\Theta})(k,t)$ is a classical solutions to (\ref{B}). 
 
\begin{Lemma}\label{2.9.} For any solutions $(\hat{H}, \hat{S})$ of (\ref{IE-p-1}) such that $||(\hat{H},\hat{S})(\cdot,p)||_{N}\in L^1(e^{-\omega p}dp)$ the Laplace transform
\begin{equation*}
(\hat{u},\hat{\Theta})(k,t)=(\hat{u}_0,\hat{\Theta}_0)(k)+\int_0^{\infty}(\hat{H},\hat{S})(k,p)e^{-p/t}dp
\end{equation*}
solves (\ref{FB2}) for $\Re (1/t)>\omega$. Moreover, $(\hat{u},\hat{\Theta})(k,t)$ is analytic for $t\in(0,\omega^{-1})$. 
\end{Lemma}
\begin{pf} Recall (\ref{Chap2H}),
\begin{equation*}
\mathcal{H}^{(\nu)}(p,p',k)=\int_{p'/p}^1\left\{\frac{1}{2\pi i}\int_{c-i\infty}^{c+i\infty}\tau^{-1}exp[-\nu|k|^2\tau^{-1}(1-s)+(p-p's^{-1})\tau]d\tau\right\}ds.
\end{equation*}
Let $\hat{G}_1=-ik_j\hat{G}_j^{[1]}+P_k(ae_2\hat{S})$ and $\hat{G}_2=-ik_j\hat{G}_j^{[2]}$. Changing variable $p'/s\rightarrow p'$ and applying Fubini's theorem gives
\begin{align}\label{long time 4.35}
\int_0^p&\left(\mathcal{H}^{(\nu)}(p,p',k)\hat{G}_1(k,p'),\mathcal{H}^{(\mu)}(p,p',k)\hat{G}_2(k,p')\right)dp'\\ \nonumber
&=\int_0^1s\left\{\int_0^p\left(\hat{G}_1(k,p's)\mathcal{I}^{(\nu)}(p-p',s,k),\hat{G}_2(k,p's)\mathcal{I}^{(\mu)}(p-p',s,k)\right)dp'\right\}ds,
\end{align}
where for $p>0$ 
\begin{equation*}
\mathcal{I}^{(\nu)}(p,s,k)=\frac{1}{2\pi i}\int_{c-i\infty}^{c+i\infty} \tau^{-1}exp[-\nu|k|^2\tau^{-1}(1-s)+p\tau]d\tau.
\end{equation*}
Taking the Laplace transform of (\ref{long time 4.35}) with respect to $p$ and again using Fubini's theorem yields
\begin{align*}
\int_0^{\infty}e^{-pt^{-1}}\int_0^1\int_0^p\left(\hat{G}_1(k,p's)\mathcal{I}^{(\nu)}(p-p',s,k),\hat{G}_2(k,p's)\mathcal{I}^{(\mu)}(p-p',s,k)\right)sdp'dsdp\\ 
=\int_0^1\left(\hat{g}_1(k,st)I^{(\nu)}(t,s,k),\hat{g}_2(k,st)I^{(\mu)}(t,s,k)\right)ds,
\end{align*}
where $\hat{g}(k,t)=\mathcal{L}[\hat{G}(k,\cdot)](t^{-1})$ and $I(t,s,k)=\mathcal{L}[\mathcal{I}(\cdot,s,k)](t^{-1})$. By assumption, $||(\hat{H},\hat{S})(\cdot,p)||_{N}\in L^1(e^{-\omega p}dp)$ and $||(\hat{u}_0,\hat{\Theta}_0)||_{N}<\infty$. From the definition of $\hat{G}^{[l]}_j$ given in (\ref{G_j}) and 
Lemma \ref{2.6.} it follows that $\hat{G}$ are Laplace transformable in p, for $t\in(0,\omega^{-1})$. Thus,
\begin{align*}
\hat{g}_1&:=-i k_j P_k[\hat{h}_j\hat{*}\hat{h}+\hat{h}_j\hat{*}\hat{u}_0+
\hat{u}_{0,j}\hat{*}\hat{h}]+P_k[ae_2\hat{s}]\\ 
\hat{g}_2&:=-i k_j[\hat{h}_j\hat{*}\hat{s}+\hat{h}_j\hat{*}\hat{\Theta}_0+
\hat{u}_{0,j}\hat{*}\hat{s}].
\end{align*}
We also have
\begin{equation*}
I^{(\nu)}(t,s,k)=te^{-\nu|k|^2t(1-s)}.
\end{equation*}
Recalling the integral equations for $(\hat{H},\hat{S})$ given in (\ref{IE-p-1}), we have  
\begin{align*}
(\hat{h},\hat{s})(k,t)-&\left(\hat{u}_1(k)\left(\frac{1-e^{-\nu|k|^2t}}{\nu|k|^2}\right),\hat{\Theta}_1(k)\left(\frac{1-e^{-\mu|k|^2t}}{\mu|k|^2}\right)\right)\\ 
&=t\int_0^{1}\left(e^{-\nu|k|^2t(1-s)}\hat{g}_1(k,st),e^{-\nu|k|^2t(1-s)}\hat{g}_2(k,st)\right)ds\\ 
&=\int_0^{t}\left(e^{-\nu|k|^2(t-s)}\hat{g}_1(k,s),e^{-\nu|k|^2(t-s)}\hat{g}_2(k,s)\right)ds. 
\end{align*}
Therefore, we directly verify $(\hat{u},\hat{\Theta})(k,t)=(\hat{u}_0,\hat{\Theta}_0)(k)+(\hat{h},\hat{s})(k,t)$ satisfies (\ref{FB2}).  Moreover, analyticity in $t$ follows from the representation
\begin{equation*} 
(\hat{u},\hat{\Theta})(k,t)=(\hat{u}_0,\hat{\Theta}_0)(k)+\int_0^{\infty}\left(\hat{H}, \hat{S}\right)(k,p)e^{-p/t}dp.
\end{equation*}
\end{pf}
\begin{Lemma}\label{2.11.} (Instantaneous smoothing) Assume $||(\hat u_0, \hat{\Theta}_0)||_{N}<\infty $ and $|| \hat f||_{N}<\infty$ with $N$ either $L^1\cap L^{\infty}(\mathbb{R}^d)$ or $(\gamma, \beta)$ with $\gamma >d$, $\beta \geq 0$. For the solution $(\hat u, \hat {\Theta})$ known to exist by Lemma \ref{2.8.} for $t\in(0,T]$ with $T<\omega^{-1}$, we have $||(1+|\cdot|)^2(\hat{u},\hat{\Theta})(\cdot, t)||_{N}<\infty$ for $t\in(0,T]$.
\end{Lemma}

\begin{pf} Our goal is to bootstrap using derivatives of $(u, \Theta)$. Consider the time interval $[\epsilon,T]$ for $\epsilon> 0$ and $T<\omega^{-1}$. Define 
\begin{equation*}
\hat{V}_\epsilon(k)=\sup_{\epsilon\leq t\leq T}|(\hat u, \hat {\Theta})|(k,t).
\end{equation*}
Since $|(\hat{u}, \hat{\Theta})(k,t)|\leq |(\hat{u}_0, \hat{\Theta}_0)(k)|+\int_0^{\infty}|(\hat{H}, \hat{S})(k,p)|e^{-\omega p}dp$,
\begin{equation*}
||\hat{V}_\epsilon(k)||_{N}\leq ||(\hat{u}_0, \hat{\Theta}_0)(k)||_{N}+||(\hat{H}, \hat{S})(k,p)||_{1}^{\omega}<\infty.
\end{equation*}
On $[\epsilon,T]$ for $\epsilon >0$, 
\begin{align*}
\hat{u}(k,t)&=e^{-\nu|k|^2t}\hat u_0(k)- \int_0^t e^{-\nu|k|^2(t-\tau)}\left(i k_j P_k[\hat u_j\hat *\hat u]+aP_k[e_2\hat{\Theta}]-\hat f\right)d\tau\\ 
\hat \Theta(k,t)&=e^{-\mu|k|^2t}\hat \Theta_0(k)-i k_j \int_0^t e^{-\mu|k|^2(t-\tau)}\left\{ (\hat u_j\hat *\hat \Theta)(k,\tau)\right\}d\tau.
\end{align*}

\noindent Therefore,
\begin{align*} 
|k||(\hat{u}, \hat{\Theta})(k, t)|&\leq \left|(\hat{u}_0, \hat{\Theta}_0)(k)\right|\sqrt{\min(\nu ,\mu) }\sup_{z\geq 0}ze^{-z^2}+|\hat{f}|\int_0^t|k|e^{-\min(\nu ,\mu)|k|^2(t-\tau)}d\tau\\
&+\left(\hat{V}_0+\hat{V}_{0}\hat{*}\hat{V}_{0}\right)\int_0^t |k|^2e^{-\min(\nu ,\mu)|k|^2(t-\tau)}d\tau.
\end{align*}
Noticing that 
\begin{equation*}
\int_0^t |k|^2e^{-\min(\nu ,\mu)|k|^2(t-\tau)}d\tau\leq \frac{1}{\min(\nu ,\mu)}
\end{equation*}
and 
\begin{equation*}
\int_0^t|k|e^{-\min(\nu ,\mu )|k|^2(t-\tau)}d\tau\leq\sup_{z\geq 0}\frac{1-e^{-z}}{\sqrt{z}}\sqrt{\frac{T}{\min(\nu ,\mu )}},
\end{equation*}
it follows that 
\begin{equation*}
\left\| |k|\hat{V}_{\epsilon/2}\right\|_{N}\leq\frac{C}{\epsilon^{1/2}}||(\hat{u}_0, \hat{\Theta}_0)||_{N}+\frac{1}{\min(\nu ,\mu)}\left(C_0||\hat{V}_0||_{N}^2+||\hat{V}_0||_{N}+C\sqrt{T}||\hat{f}||_{N}\right)<\infty.
\end{equation*}
In the same spirit, for $t\in[\frac{\epsilon}{2},T]$, we have
\begin{align*} 
\hat{u}(k,t)&=e^{-\nu|k|^2t}\hat u(k,\epsilon/2)- \int_{\epsilon/2}^t e^{-\nu|k|^2(t-\tau)}\left(P_k(\hat u_j\hat *[ik_j\hat u]+ae_2\hat{\Theta})(k,\tau)-\hat f(k)\right)d\tau\\
\hat \Theta(k,t)&=e^{-\mu|k|^2t}\hat \Theta(k,\epsilon/2)-i \int_{\epsilon/2}^t e^{\frac{-|k|^2(t-\tau)}{\mu \sigma}}\left\{ (\hat u_j\hat *k_j\hat{\Theta})(k,\tau)\right\}d\tau,
\end{align*}
where we used the divergence free conditions $k\cdot\hat{u}=0$. Multiplying by $|k|^2$ and using our previous bounds, we have for $t\in[\epsilon,T]$
\begin{align*}
|k|^2|(\hat{u}, \hat{\Theta})(k, t)|\leq &\left|(\hat{u}, \hat{\Theta})(k,\epsilon/2)\right|\frac{1}{(t-\epsilon/2)\min(\nu,\mu)}\sup_{z\geq0}ze^{-z}\\ 
&+(\hat{V}_{\epsilon/2}\hat{*}|k|\hat{V}_{\epsilon/2}+|k|\hat{V}_{\epsilon/2}+|\hat{f}|) \int_{\epsilon/2}^t |k|^2e^{-\min(\nu ,\frac{1}{\mu \sigma})|k|^2(t-\tau)}d\tau
\end{align*}
Hence, 
\begin{equation*}
\left\| |k|^2\hat{V}_{\epsilon}\right\|_{N}\leq\frac{C}{\epsilon}||(\hat{u}_0, \hat{\Theta}_0)||_{N}+\frac{\left(\left(C_0\left\|\hat{V}_{\epsilon/2}\right\|_{N}+1\right)\left\||k|\hat{V}_{\epsilon/2}\right\|_{N}+||\hat{f}||_{N}\right)}{ \min(\nu ,\mu)}.
\end{equation*}
All the terms on the right hand side are bounded, which gives $||(1+|k|)^2\hat{V}_{\epsilon}||_{N}<\infty$. Further, as $\epsilon>0$ is arbitrary, it follows that $||(1+|\cdot|)^2(\hat{u}, \hat{\Theta})(\cdot,t)||_{N}<\infty$ for $t\in(0,T]$.
\end{pf}
\begin{Remark}\label{remsmooth}
{\rm 
We note that the smoothness argument in $x$ 
of the previous Lemma
can be easily extended further 
to show $\left \| (1+|k|)^4 {\hat V}_{\epsilon} 
\right \|_N$ is finite provided $\| (1+|k|^2) {\hat f} \|_N$,
is finite. Since
$\epsilon > 0$ is arbitrary, this implies instantaneous smoothing
two orders more than the forcing.
}
\end{Remark}  

\begin{Lemma}\label{equivalence}
Given $(\hat{u},\hat{\Theta})$ a solution to (\ref{FB2}) such that $||(1+|\cdot|)^2(\hat{u},\hat{\Theta})(\cdot, t)||_N<\infty$ for $t\in (0,\omega^{-1})$, then $(u,\Theta)\in L^{\infty}[0,\omega^{-1},H^2(\mathbb{R}^d)]$ solves (\ref{B}). 
\end{Lemma}

\begin{pf}
Suppose $(\hat{u},\hat{\Theta})$ is a solution to (\ref{FB2}) such that $||(1+|\cdot|)^2(\hat{u},\hat{\Theta})(\cdot, t)||_N<\infty$ for $t\in (0,\omega^{-1})$. We notice that by our choice of norms, $(1+|\cdot|)^2(\hat{u},\hat{\Theta})(\cdot,t)\in L^2(\mathbb{R}^d)$ for any $t\in (0,\omega^{-1})$. Indeed for $N=(\gamma, \beta)$, we have
\begin{equation*}
\int(1+|k|)^4|(\hat{u},\hat{\Theta})(k,t)|^2dk\leq
||(1+|\cdot|)^2|(\hat{u},\hat{\Theta})(\cdot,t)||_{\gamma, \beta}^2\int\frac{e^{-2\beta|k|}}{(1+|k|)^{2\gamma}}dk.
\end{equation*}
As $\gamma>d$, $\int\frac{1}{(1+|k|)^{2\gamma}}e^{-2\beta|k|}dk<\infty$. For $N=L^1\cap L^{\infty}$ we have,
\begin{equation*}
\int(1+|k|)^4|(\hat{u},\hat{\Theta})(k,t)|^2dk\leq\int(1+|k|)^2|(\hat{u},\hat{\Theta})(k,t)|dk\sup_{k\in \mathbb{R}^d}(1+|k|)^2|(\hat{u},\hat{\Theta})(k,t)|.
\end{equation*}
So, $||(1+|\cdot|)^2(\hat{u},\hat{\Theta})(\cdot,t)||_{L^2(\mathbb{R}^d)}\leq ||(1+|\cdot|)^2(\hat{u},\hat{\Theta})(\cdot,t)||_{L^1\cap L^{\infty}(\mathbb{R}^d)}$. Thus, by well known properties of the Fourier transform $(u,\Theta)=\mathcal{F}^{-1}(\hat{u},\hat{\Theta})(x,t)\in L^{\infty}(0,\omega^{-1},H^2(\mathbb{R}^d))$. As $(\hat{u},\hat{\Theta})$ solves (\ref{FB2}), $(\hat{u},\hat{\Theta})$ is differentiable almost everywhere and 
\begin{align*}
\hat{u}_t+\nu |k|^2\hat{u}&=-i k_j P_k[\hat{u}_j\hat{*}\hat{u}]+a P_k[e_2\hat{\Theta}]+\hat{f}\\ 
\hat{\Theta}_t+\mu |k|^2\hat{\Theta}&=-i k_j[\hat{u}_j\hat{*}\hat{\Theta}], \quad k\in \mathbb{R}^d \quad t\in \mathbb{R}^{+}.
\end{align*}
Further, $(\hat{u}_t,\hat{\Theta}_t)(k,t)\in L^{\infty}(0,\omega^{-1},L^2(\mathbb{R}^d))$ since $(1+|k|)^2(\hat{u},\hat{\Theta})(k,t)\in L^{\infty}(0,\omega^{-1},L^2(\mathbb{R}^d))$. Hence, $(u,\Theta)(x,t)=\mathcal{F}^{-1}(\hat{u},\hat{\Theta})(x,t)$ solves 
\begin{align*}
u_t-\nu\Delta u&=-P[ u\cdot \nabla u - a e_2\Theta] +f(x)\\ 
\Theta_t-\mu \Delta \Theta&= -u\cdot \nabla \Theta.
\end{align*}
\end{pf}

\noindent {\bf Proof of Theorem \ref{existence}:} 
Suppose $||(1+|\cdot|)^2(\hat{u}_0, \hat{\Theta}_0)||_{N}<\infty$ and $||\hat{f}||_{N}<\infty$. Then from the definition of $(\hat u_1, \hat \Theta_1)$ in $(\ref{u1})$ we see $||(\hat{u}_1,\hat{\Theta}_1)||_{N}<\infty$, since
\begin{align*}
||(\hat u_1,\hat \Theta_1)||_{N}\leq \max(\nu, \mu) \left\||k|^2(\hat u_0, \hat{\Theta}_0)\right\|_{N}&+ C_0||\hat u_0||_{N}\left\||k|(\hat u_0,\hat \Theta_0)\right\|_{N}\\
&+ a||\hat \Theta_0||_{N}+||\hat f||_{N}.
\end{align*}
Therefore, when $\omega$ is large enough to ensures ($\ref{1-2.39}$), Lemma \ref{2.8.} gives $(\hat H, \hat S)(k, \cdot)$ is in $L^1(e^{-\omega p}dp)$. Applying Lemma \ref{2.9.}, we know for $t$ such that $\Re \frac{1}{t}>\omega$, $(\hat{H},\hat{S})(k,p)$ is Laplace transformable in $1/t$ with $(\hat{u},\hat{\Theta})(k,t)=(\hat{u}_0,\hat{\Theta}_0)(k)+(\hat{h},\hat{s})(k,t)$
satisfying Boussinesq equation in the Fourier space, (\ref{FB2}). Since $||(\hat H, \hat S)(\cdot, p)||_{N}<\infty$, we have $||(\hat u, \hat \Theta)(\cdot, t)||_{N}<\infty$ if $\Re \frac{1}{t}>\omega$, and i) is proved. Moreover, Lemma \ref{2.9.} shows that $(\hat{u},\hat{\Theta})$ is analytic for $\Re \frac{1}{t}>\omega$ and has the representation
\begin{equation*}
(\hat{u}, \hat{\Theta})(k,t)=(\hat{u}_0, \hat{\Theta}_0)(x)+\int_0^{\infty}(\hat{H}, \hat{S})(k,p)e^{-p/t}dp
\end{equation*}
proving ii). For iii), Lemma \ref{2.11.} shows that $||(1+|\cdot|)^2(\hat{u}, \hat{\Theta})(\cdot,t)||_{N} < \infty$ for $t\in[0,\omega^{-1})$ while Lemma \ref{equivalence} shows that $(u, \Theta)(x,t)\in L^{\infty}(0,T,H^2(\mathbb{R}^d))$ solves (\ref{B}). Moreover, $(u, \Theta)(x,t)$ is the unique solution to (\ref{B}) in $L^{\infty}(0,T,H^2(\mathbb{R}^d))$ as classical solutions are known to be unique, \cite{Temam}. Finally, suppose $(\hat H, \hat S)(k, \cdot)$ is in $L^1(e^{-\omega p}dp)$ for any $\omega >0$. By Lemma \ref{2.9.}, we know for any $t>0$, $(\hat{H},\hat{S})(k,p)$ is Laplace transformable with $(\hat{u},\hat{\Theta})(k,t)=(\hat{u}_0,\hat{\Theta}_0)(k)+(\hat{h},\hat{s})(k,t)$
satisfying Boussinesq equation in the Fourier space, (\ref{FB2}). Further, appealing to instantaneous smoothing Lemma \ref{2.11.} the solution is smooth. Thus, if $(\hat H, \hat S)(k, \cdot)$ is in $L^1(e^{-\omega p}dp)$ for any $\omega >0$, then a smooth global solution exists and iv) is proved. 

\section{Borel-Summability}\label{sec6}

\indent We now show Borel-summability of the solutions guaranteed by Theorem \ref{existence} for $\beta >0$. This requires us to show that the solutions $(\hat{H}, \hat{S})(k,p)$ to the Boussinesq equation in Borel space is analytic in $p$ for $p\in \left\{0\right\}\cup\mathbb{R}^+$. First, we will seek a solution which is a power series

\begin{equation}\label{1-3.49}
(\hat{H}, \hat{S})(k,p)-(\hat{u}_1,\hat{\Theta}_1)(k)=\sum_{l=1}^{\infty}(\hat{H}^{[l]}, \hat{S}^{[l]})(k)p^l.
\end{equation} 

\begin{Remark}
{\rm 
We will use induction to bound the successive terms of the power series. Many of these bounds have constants depending on the dimension in $k$ as before. For brevity of notation the dependence on dimension is suppressed after introducing the constants.
}
\end{Remark}

For the purpose of finding power series solutions, (\ref{1-2.18}) is not a good representation. By construction, $\frac{\pi}{z} \mathcal{G}(z,z')$ satisfies $[p\partial_{pp}+2\partial_p+\nu|k|^2]y=0$ with $\frac{\pi}{z} \mathcal{G}(z,z')\rightarrow 0$ and $\partial_p\left(\frac{\pi}{z} \mathcal{G}(z,z')\right)\rightarrow \frac{1}{p}$ as $p'$ approaches $p$ from below. Hence, we have the equivalent equations
\begin{align}\label{1-2.9}
[p\partial_{pp}+2\partial_p+\nu|k|^2]\hat{H}&=ik_j\hat{G}_j^{[1]}+aP_k[\hat{e}_2\hat{S}]\\ \nonumber
[p\partial_{pp}+2\partial_p+\mu|k|^2]\hat{S}&=ik_j\hat{G}_j^{[2]}.
\end{align}
\noindent We substitute ($\ref{1-3.49}$) into ($\ref{1-2.9}$) and identify powers of $p^l$ to get a relationship for the coefficients. We will use the fact that
\begin{equation*}
p^l*p^n=\frac{l!n!}{(l+n+1)!}p^{l+n+1}.
\end{equation*}
For $l=0$, we have
\begin{align}\label{1-3.50}
2\hat{H}^{[1]}&=-i k_j P_k[\hat{u}_{1,j}\hat{*}\hat{u}_0+\hat{u}_{0,j}\hat{*}\hat{u}_1]-\nu |k|^2\hat{u}_1+P_k[a e_2\hat{\Theta}_1]\\ \nonumber
2\hat{S}^{[1]}&=-i k_j[\hat{u}_{1,j}\hat{*}\hat{\Theta}_0+\hat{u}_{0,j}\hat{*}\hat{\Theta}_1]-\mu |k|^2\hat{\Theta}_1.
\end{align}
For $l=1$, we have
\begin{align}\label{1-3.51}
6\hat{H}^{[2]}+\nu |k|^2\hat{H}^{[1]}&=-i k_j P_k[\hat{H}^{[1]}_j\hat{*}\hat{u}_0+\hat{u}_{0,j}\hat{*}\hat{H}^{[1]}+\hat{u}_{1,j}\hat{*}\hat{u}_1]+P_k[a e_2\hat{S}^{[1]}] \\ \nonumber
6\hat{S}^{[2]}+\mu |k|^2\hat{S}^{[1]}&=-i k_j[\hat{S}^{[1]}_j\hat{*}\hat{\Theta}_0+\hat{u}_{0,j}\hat{*}\hat{S}^{[1]}+\hat{u}_{1,j}\hat{*}\hat{\Theta}_1].
\end{align}
\noindent More generally, for $l\geq 2$, we have
\begin{align}\label{1-3.53}
(l+1)(l+2)\hat{H}^{[l+1]}&=-\nu|k|^2\hat{H}^{[l]}-i k_j P_k\left[\sum_{l_1=1}^{l-2}\frac{l_1!(l-l_1-1)!}{l!}\hat{H}_j^{[l_1]}\hat{*}\hat{H}^{[l-l_1-1]}\right]\\ \nonumber 
-i k_j P_k[\hat{u}_{0,j}&\hat{*}\hat{H}^{[l]}+ \hat{H}_j^{[l]}\hat{*}\hat{u}_0+\frac{1}{l}\hat{u}_{1,j}\hat{*}\hat{H}^{[l-1]} +\frac{1}{l}\hat{H}_j^{[l-1]}\hat{*}\hat{u}_1] +P_k[a e_2\hat{S}^{[l]}]\\ \label{1-3.53b}
(l+1)(l+2)\hat{S}^{[l+1]}&=-\mu|k|^2\hat{S}^{[l]}-i k_j\left[\sum_{l_1=1}^{l-2}\frac{l_1!(l-l_1-1)!}{l!}\hat{H}_j^{[l_1]}\hat{*}\hat{S}^{[l-l_1-1]}\right]\\ \nonumber
-i k_j\left[\right.\hat{u}_{0,j}&\hat{*}\hat{S}^{[l]}+\hat{H}_j^{[l]}\hat{*}\hat{\Theta}_0+ \frac{1}{l}\hat{u}_{1,j}\hat{*}\hat{S}^{[l-1]}+ \frac{1}{l}\hat{H}_j^{[l-1]}\hat{*}\hat{\Theta}_1\left. \right].
\end{align}

\begin{Definition} It is useful to define a n-th order polynomial, call it $\mathcal{Q}_n$, 
\begin{equation*}
\mathcal{Q}_n(y)=\sum_{j=0}^{n}2^{n-j}\frac{y^j}{j!}.
\end{equation*}
\end{Definition}
\begin{Definition}\label{maximums}
It is also useful to define the constant
\begin{equation*}
M_1=\max(\nu,\mu).
\end{equation*}
\end{Definition}

\subsection{Estimates on the Solution in the Borel Plane}
\begin{Lemma}\label{3.2.} If $||(\hat{u}_0, \hat{\Theta}_0)||_{\gamma+2,\beta}<\infty$ for $\gamma >d$ and $\beta>0$, then there are constants $A_0$, $D_0>0$ not depending on $l$ or $k$ such that
\begin{equation}\label{1-3.54}
|(\hat{H}^{[l]}, \hat{S}^{[l]})|\leq e^{-\beta |k|}A_0D_0^l(1+|k|)^{-\gamma}\frac{\mathcal{Q}_{2l}(|\beta k|)}{(2l+1)^2}.
\end{equation}
Furthermore, the solutions guaranteed to exist in Lemma (\ref{2.8.}) have convergent power series representations in $p$, and for $|p|<(4D_0)^{-1}$ 
\begin{equation*}
(\hat{H}, \hat{S})(k,p)=(\hat{u}_1, \hat{\Theta}_1)(k)+\sum_{l=1}^{\infty}(\hat{H}^{[l]}, \hat{S}^{[l]})(k)p^l. \end{equation*}
\end{Lemma}

To prove this lemma we will establish bounds for $(\hat{H}^{[l]},\hat{S}^{[l]})$ using induction. 

\begin{Lemma}\label{3.3.} For the base case, we have
\begin{equation}\label{1-3.55}
|(\hat{H}^{[1]}, \hat{S}^{[1]})(k)|\leq \frac{e^{-\beta |k|}\mathcal{Q}_{2}(\beta |k|)A_0D_0}{(1+|k|)^{\gamma}9}
\end{equation}
\noindent for
\begin{equation*}
A_0D_0\geq \frac{9}{\beta^2}||(\hat{u}_1, \hat{\Theta}_1)||_{\gamma,\beta}\left(C_0\beta||(\hat{u}_0, \hat{\Theta}_0)||_{\gamma,\beta}+ M_1+a\beta^2\right)
\end{equation*}
\end{Lemma}

\begin{pf} From $(\ref{1-3.50})$ and Lemma \ref{2.3.}, we get
\begin{align}\label{1-3.52}
|(\hat{H}^{[1]}, \hat{S}^{[1]})(k)|\leq &\frac{e^{-\beta|k|}}{2(1+|k|)^{\gamma}}\left( |k|^2||(\hat{u}_1, \hat{\Theta}_1)||_{\gamma, \beta}M_1\right. \\ \nonumber
 &+\left. 2C_0|k|\left\|(\hat{u}_0, \hat{\Theta}_0)\right\|_{\gamma, \beta}||(\hat{u}_1, \hat{\Theta}_1)||_{\gamma, \beta}+ a||\hat{\Theta}_1||_{\gamma, \beta}\right).
\end{align}
The result now follows after noting that $\mathcal{Q}_2(\beta |k|)=4+2\beta |k|+\frac{1}{2}(\beta |k|)^2$.
\end{pf}

\noindent For the general terms we will need a series of lemmas, which depend on the Fourier inequalities developed in \ref{sec8}, bounding the terms that appear on the right side of (\ref{1-3.53}). 

\begin{Lemma}\label{3.4.} Assume that $(\hat{H}^{[l]}, \hat{S}^{[l]})$ satisfies (\ref{1-3.54}) for $l\geq 1$. Then we have,
\begin{equation*}
\frac{|k|^2|(\hat{H}^{[l]}, \hat{S}^{[l]})|}{(l+1)(l+2)}\leq\frac{6A_0D_0^l e^{-\beta|k|}\mathcal{Q}_{2l+2}(\beta|k|)}{\beta^2(1+|k|)^{\gamma}(2l+3)^2}.
\end{equation*}
\end{Lemma}
\begin{pf} The proof follows from (\ref{1-3.54}) directly by noting that for $y\geq 0$
\begin{equation*}
\frac{y^2\mathcal{Q}_{2l}(y)}{(2l+2)(2l+1)}\leq \mathcal{Q}_{2l+2}(y) \textnormal{ and }\frac{(2l+2)(2l+3)^2}{(l+1)(l+2)(2l+1)}\leq 6.
\end{equation*}
\end{pf}
\begin{Lemma}\label{3.5.} Suppose $(\hat{H}^{[l]}, \hat{S}^{[l]})$ satisfies (\ref{1-3.54}) for $l\geq 1$. Then both
\begin{equation*}
\left|k_j\left(P_k(\hat{u}_{0,j}\hat{*}\hat{H}^{[l]}), \hat{u}_{0,j}\hat{*}\hat{S}^{[l]}\right)\right| \textnormal{ and }
\left|k_j\left(P_k(\hat{H}^{[l]}_j\hat{*}\hat{u}_{0}), \hat{H}^{[l]}_j\hat{*}\hat{\Theta}_0\right)\right|
\end{equation*}
are bounded by 
\begin{equation*}
 2^{\gamma}||(\hat{u}_0, \hat{\Theta}_0)||_{\gamma,\beta}\frac{9C_7\pi A_0D_0^l e^{-\beta|k|}(l+1)(l+2)}{2\beta^d(2l+3)^2(1+|k|)^\gamma}\mathcal{Q}_{2l+2}(|\beta k|).
\end{equation*}
Similarly, suppose $(\hat{H}^{[l-1]},\hat{S}^{[l-1]})$ satisfies (\ref{1-3.54}) for $l\geq 2$. Then both
\begin{equation*}
\left|k_j\left(P_k(\hat{u}_{1,j}\hat{*}\hat{H}^{[l-1]}), \hat{u}_{1,j}\hat{*}\hat{S}^{[l-1]}\right)\right|\textnormal{ and } \left|k_j\left(P_k(\hat{H}^{[l-1]}_j\hat{*}\hat{u}_{1}),\hat{H}^{[l-1]}_j\hat{*}\hat{\Theta}_1\right)\right|
\end{equation*}
are bounded by
\begin{equation*}
2^{\gamma}||(\hat{u}_1, \hat{\Theta}_1)||_{\gamma,\beta}\frac{9C_7\pi A_0 D_0^{l-1} e^{-\beta|k|}l(l+1)\mathcal{Q}_{2l}(|\beta k|)}{2\beta^d(2l+1)^2(1+|k|)^\gamma}.
\end{equation*}
\end{Lemma}
\begin{pf} We use the estimate (\ref{1-3.54}) on $(\hat{H}^{[l]}, \hat{S}^{[l]})$ and \ref{6.7.} in $\mathbb{R}^d$ with $n=0$ to get
\begin{align*}
|k_j\hat{u}_{0,j}\hat{*}(\hat{H}^{[l]},\hat{S}^{[l]})|&\leq ||\hat{u}_0||_{\gamma, \beta}\frac{A_0D_0^l}{(2l+1)^2}\left(|k|\int_{k'\in\mathbb{R}^d}\frac{e^{-\beta(|k'|+|k-k'|)}}{(1+|k'|)^{\gamma}(1+|k-k'|)^{\gamma}}\mathcal{Q}_{2l}(\beta|k'|)dk'\right)\\ 
&\leq \frac{||\hat{u}_0||_{\gamma, \beta}A_0D_0^l}{(2l+1)^2}\sum_{m=0}^{2l}\frac{2^{2l-m}}{m!}|k|\int_{k'\in\mathbb{R}^d}\frac{e^{-\beta(|k'|+|k-k'|)}}{(1+|k'|)^{\gamma}(1+|k-k'|)^{\gamma}}|\beta k'|^m dk'\\  
&\leq \frac{C_7\pi ||\hat{u}_0||_{\gamma, \beta}A_0D_0^l2^{\gamma}e^{-\beta|k|}}{(2l+1)^2\beta^d(1+|k|)^{\gamma}}\sum_{m=0}^{2l}2^{2l-m}(m+2)\mathcal{Q}_{m+2}(\beta|k|)\\ 
&\leq \frac{2^{\gamma}C_7\pi ||\hat{u}_0||_{\gamma, \beta}A_0D_0^l e^{-\beta|k|}}{(2l+1)\beta^d(1+|k|)^{\gamma}}(l+2)\mathcal{Q}_{2l+2}(\beta|k|).
\end{align*}
\noindent The first part of the lemma now follows noting $\frac{2(2l+3)^2}{(2l+1)(l+1)}\leq 9$ for $l\geq 1$. The second parts is proved similarly.
\end{pf}

\begin{Lemma}\label{3.7.} Let $l\geq 3$. Suppose $(\hat{H}^{[l_1]}, \hat{S}^{[l_1]})$ and $(\hat{H}^{[l-1-l_1]}, \hat{S}^{[l-1-l_1]})$ satisfy $(\ref{1-3.54})$ for $l_1=1,\dots, l-2$. Then
\begin{equation*}
\left|k_j\left[\sum_{l_1=1}^{l-2}\frac{l_1!(l-1-l_1)!}{(l+2)!}\left(P_k(\hat{H}_j^{[l_1]}\hat{*}\hat{H}^{[l-1-l_1]}),\hat{H}_j^{[l_1]}\hat{*}\hat{S}^{[l-1-l_1]}\right)\right]\right| 
\end{equation*}
is bounded by
\begin{equation*}
2^{\gamma+3} C_7 A_0^2D_0^{l-1}(1+|k|)^{-\gamma}e^{-\beta|k|}\frac{\mathcal{Q}_{2l}(\beta|k|)}{\beta^d (2l+3)^2}.
\end{equation*} 
\end{Lemma}

\begin{pf} The proof is similar to that in \cite{smalltime} with $\hat{W}^{[l_2]}$ replaced by $(\hat{W}^{[l_2]}, \hat{Q}^{[l_2]})$. For more details see \cite{smalltime} and \cite{Thesis}.
\end{pf}

\begin{Lemma}\label{3.8.} For $l=2$ we have,
\begin{align*}
|(\hat{H}^{[2]},\hat{S}^{[2]})|\leq& \frac{e^{-\beta|k|}\mathcal{Q}_4(\beta |k|)}{5^2(1+|k|)^{\gamma }}\left(\frac{6A_0D_0M_1}{\beta^2}+\frac{2^{\gamma}9C_7\pi A_0D_0||(\hat{u}_0, \hat{\Theta}_0)||_{\gamma,\beta}}{\beta^d}\right.\\ 
&\left. \quad+A_0D_0a+\frac{C_0}{\beta}||(\hat{u}_1, \hat{\Theta}_1)||^2_{\gamma, \beta}\right).
\end{align*}
\noindent Thus, $(\hat{H}^{[2]}, \hat{S}^{[2]})$ satisfies ($\ref{1-3.54}$) for 
\begin{equation}\label{1-3.57}
D_0^2\geq\frac{6D_0M_1}{\beta^2}+D_0a+\frac{2^{\gamma}9C_7\pi D_0}{\beta^d}||(\hat{u}_0, \hat{\Theta}_0)||_{\gamma,\beta}+ \frac{C_0}{A_0\beta}||(\hat{u}_1, \hat{\Theta}_1)||^2_{\gamma, \beta}.
\end{equation}
\end{Lemma}
\begin{pf} We start from $(\ref{1-3.51})$. For the first term we use Lemma \ref{3.4.}. For the second term, appearing in (\ref{1-3.57}), we use our induction assumption and $\frac{Q_2(\beta |k|)}{54}\leq\frac{Q_4(\beta|k|)}{25}$. For the next term, we use Lemma \ref{3.5.}. For the last terms, apply Corollary \ref{2.2.} and use $\frac{|k|}{6}\leq\frac{Q_4(\beta|k|)}{25 \beta}$.
\end{pf}

\noindent{\bf Proof of Lemma \ref{3.2.}}
The base case is proved picking $D_0$ large enough so $(\ref{1-3.57})$ and $(\ref{1-3.55})$ hold. For general $l\geq 2$ suppose $(\hat{H}^{[m]}, \hat{S}^{[m]})$ satisfies (\ref{1-3.54}) for $m=1, \dots, l$. We estimate terms on the right of $(\ref{1-3.53})$ and $(\ref{1-3.53b})$, using Lemma \ref{3.4.}, \ref{3.5.}, and \ref{3.7.} and the fact that $Q_{2l}(y)\leq 1/4Q_{2l+2}(y)$, to get 
\begin{align*}
|(\hat{H}^{[l+1]},\hat{S}^{[l+1]})|&\leq \frac{A_0D_0^{l-1}Q_{2l+2}(\beta |k|)}{(2l+3)^2(1+|k|)^{\gamma}}\left\{\frac{6D_0M_1}{\beta^2}+\frac{a D_0}{2}+\frac{2^{\gamma}9C_7\pi D_0}{\beta^d}||(\hat{u}_0, \hat{\Theta}_0)||_{\gamma, \beta} \right.  \\ 
 &\left. \qquad +\frac{2^{\gamma}9C_7\pi (2l+3)^2}{4(l+2)(2l+1)^2\beta^d}||(\hat{u}_1, \hat{\Theta}_1)||_{\gamma, \beta}+\frac{2^{\gamma +3}C_7A_0}{4\beta^d}\right\}\\ 
&\leq \frac{A_0D_0^{l+1}e^{-\beta |k|}}{(1+|k|)^{\gamma}(2l+3)^2}\mathcal{Q}_{2l+2}(\beta |k|)
\end{align*}
\noindent where $D_0$ has been chosen large enough so
\begin{align*}
\left\{  \frac{6D_0M_1}{\beta^2}\right.+\frac{a D_0}{2} +&\frac{2^{\gamma}9C_7\pi D_0}{\beta^d}||(\hat{u}_0, \hat{\Theta}_0)||_{\gamma, \beta}+\frac{2^{\gamma}9C_7\pi D_0}{4\beta^d}||(\hat{u}_1, \hat{\Theta}_1)||_{\gamma, \beta}\\ 
&\left. +\frac{2^{\gamma +1}C_7A_0}{\beta^d}\right\} \leq D_0^2.
\end{align*}
\noindent We also used $\frac{(2l+3)^2}{(2l+1)^2(l+2)}\leq 1$ in the above. Thus, by induction, we have $(\ref{1-3.54})$ satisfied for any $l\geq 1$. So, $\sum_{l=1}^{\infty}(\hat{H}^{[l]},\hat{S}^{[l]})(k)p^l$ is convergent for $|p|\leq \frac{1}{4 D_0}$ since $\mathcal{Q}_{2l}(\beta |k|)\leq 4^l e^{\beta|k|/2}$. By construction of the iteration, $(\hat{H}, \hat{S})-(\hat{u}_1,\hat{\Theta}_1)= \sum_{l=1}^{\infty}(\hat{H}^{[l]}, \hat{S}^{[l]})(k)p^l$ is a solution to $(\ref{1-2.9})$ which is zero at $p=0$.  However, we know there is a unique solutions to $(\ref{1-2.9})$ which is zero and $p=0$ in the space $\mathcal{A}_L^{\infty}$, which includes analytic functions at the origin for $L$ sufficiently small. Thus, for $(\hat{H}, \hat{S})$ the solution guaranteed by Lemma \ref{2.8.}, we have
\begin{equation*}
(\hat{H}, \hat{S})(k,p)=(\hat{u}_1, \hat{\Theta}_1)(k)+\sum_{l=1}^{\infty}(\hat{H}^{[l]}, \hat{S}^{[l]})(k)p^l.
\end{equation*}

\subsection{Estimates on $\partial_p^l(\hat{H}, \hat{S})(k,p)$}

We now want to develop estimates on $\partial_p^l(\hat{H}, \hat{S})(k,p)$ in order to show that we can analytically extend our solutions along $\mathbb{R}^{+}$ with a radius of convergence independent of center $p_0$ along $\mathbb{R}^{+}$. Combining this with the fact that the solutions are exponentially bounded will give Borel summability.
\begin{Definition}\label{Def4.1}For $l\geq 1$ we define,
\begin{align*}
(\hat{H}^{[l]}, \hat{S}^{[l]})(k,p)&=\frac{1}{l!}\partial_p^l(\hat{H}, \hat{S})(k,p) \\  (\hat{H}^{[0]}, \hat{S}^{[0]})(k,p)&=(\hat{H}, \hat{S})(k,p)-(\hat{u}_1,\hat{\Theta}_1).
\end{align*}
\end{Definition}

\begin{Lemma}\label{4.2.} If $||(\hat{u}_0, \hat{\Theta}_0)||_{\gamma+2,\beta}<\infty$ for and $\beta>0$, then there are constants $A$, $D>0$ not depending on $l, k$ or $p$ such that
\begin{equation}\label{1-4.60}
|(\hat{H}^{[l]}, \hat{S}^{[l]})(k,p)|\leq \frac{e^{\omega' p}e^{-\beta |k|}AD^l}{(1+p^2)(1+|k|)^{\gamma}}\frac{\mathcal{Q}_{2l}(|\beta k|)}{(2l+1)^2}
\end{equation}
\noindent where $\omega'=\omega+1$ for $\omega$ chosen as in Lemma \ref{2.8.}. We will prove the lemma by induction, and as before we will develop several lemmas to establish the bound. 
\end{Lemma}

For $l=0$, we use Lemma \ref{2.8.} which says that for $\omega$ sufficiently large
\begin{equation*}
|(\hat H, \hat S)(k,p)|\leq \frac{2e^{-\beta |k|+\omega p}||(\hat u_1, \hat \Theta_1)||_{\gamma, \beta}}{(1+|k|)^{\gamma}}.
\end{equation*}
\noindent We chose $\omega'=\omega+1$ and recall Definition \ref{Def4.1} to get
\begin{equation}\label{4.59}
|(\hat H^{[0]}, \hat{S}^{[0]})(k,p)|\leq \frac{3e^{-\beta |k|+\omega' p}||(\hat u_1, \hat \Theta_1)||_{\gamma, \beta}}{(1+p^2)(1+|k|)^{\gamma}},
\end{equation}
\noindent and the base cases of $(\ref{1-4.60})$ is proved for $A=3||(\hat{u}_1, \hat{\Theta}_1)||_{\gamma, \beta}$. 

For the general case ($l\geq1$) we take $\partial _p^l$ in $(\ref{1-2.9})$ and divide by $l!$, to obtain
\begin{multline}\label{1-4.62}
p\hat{H}^{[l]}_{pp}+(l+2)\hat{H}^{[l]}_p+\nu |k|^2\hat{H}^{[l]}=\left(-i k_j P_k[\hat{u}_{0,j}\hat{*}\hat{u}_1+\hat{u}_{1,j}\hat{*}\hat{u}_0]-\nu|k|^2\hat{u}_1\right)\delta_{l,0}\hspace{.5 in}\\ 
-i k_j P_k\left[ \int_0^p\hat{H}^{[l]}_j(\cdot,p-s)\hat{*}\hat{H}^{[0]}(\cdot,s)ds +\sum_{l_1=1}^{l-1}\frac{l_1!(l-l_1-1)!}{l!}\hat{H}^{[l_1]}_j(\cdot,0)\hat{*}\hat{H}^{[l-l_1-1]}(\cdot,p)\right]\\ 
-i k_j P_k[\frac{1}{l}(\hat{u}_{1,j}\hat{*}\hat{H}^{[l-1]}+\hat{H}_j^{[l-1]}\hat{*}\hat{u}_1)+\hat{H}^{[l]}_j\hat{*}\hat{u}_0+
\hat{u}_{0,j}\hat{*}\hat{H}^{[l]} + \delta_{l,1}\hat{u}_{1,j}\hat{*}\hat{u}_1]+P_k(a e_2\hat S^{[l]})
\end{multline}
\begin{multline}\label{2-4.62}
p\hat{S}^{[l]}_{pp}+(l+2)\hat S^{[l]}_p+\mu |k|^2\hat{S}^{[l]}=\left(-i k_j[\hat{u}_{0,j}\hat{*}\hat{\Theta}_1+\hat{u}_{1,j}\hat{*}\hat{\Theta}_0]-\mu|k|^2\hat{\Theta}_1\right)\delta_{l,0}\hspace{.5in}\\ 
-i k_j\left[ \int_0^p\hat{H}^{[l]}_j(\cdot,p-s)\hat{*}\hat{S}^{[0]}(\cdot,s)ds +\sum_{l_1=1}^{l-1}\frac{l_1!(l-l_1-1)!}{l!}\hat{H}^{[l_1]}_j(\cdot,0)\hat{*}\hat{S}^{[l-l_1-1]}(\cdot,p)\right]\\ -i k_j[\frac{1}{l}(u_{1,j}\hat{*}\hat{S}^{[l-1]}+\hat{H}_{j}^{[l-1]}\hat{*}\hat{\Theta}_1)+\hat{H}_j^{[l]}\hat{*}\hat{\Theta}_0+\hat{u}_{0,j}\hat{*}\hat{S}^{[l]}+ \delta_{l=1}\hat{u}_{1,j}\hat{*}\hat{\Theta}_1].
\end{multline}
\noindent Denote the right hand side of these four equations by $R_m^{[l]}$ for $m=1$ and $2$ respectively.

\begin{Lemma}\label{4.4.} For any $l\geq 0$ and for some absolute constant $C_6$, if $(\hat{H}^{[l]}, \hat{S}^{[l]})$ satisfies $(\ref{1-4.60})$, and is bounded at $p=0$ then
\begin{equation*}
|(\hat{H}^{[l+1]}, \hat{S}^{[l+1]})(k,p)|\leq \frac{C_6}{(l+1)^{5/3}}\sup_{p'\in [0,p]}|(\hat{R}^{[l]}_1, \hat{R}_2^{[l]})|+\frac{M_1|k|^2|(\hat{H}^{[l]}, \hat{S}^{[l]})(k,0)|}{(l+1)(l+2)}.
\end{equation*}
\end{Lemma}
\begin{pf} The proof is in \cite{smalltime} under Lemma 4.4. The lemma is dependent only on the operator $\mathcal{D}$ which is the same in our case. 
\end{pf}

\begin{Lemma}\label{4.6.} Suppose $(\hat{H}^{[l]}, \hat{S}^{[l]})$ satisfies $(\ref{1-4.60})$ for $l\geq 1$. Then
\begin{equation*}
\left|k_j\left(P_k(\hat{u}_{0,j}\hat{*}\hat{H}^{[l]}), \hat{u}_{0,j}\hat{*}\hat{S}^{[l]}\right)\right|\textnormal{ and }
\left|k_j\left(P_k(\hat{H}^{[l]}_j\hat{*}\hat{u}_{0}), \hat{H}^{[l]}_j\hat{*}\hat{\Theta}_0\right)\right|
\end{equation*} 
are bounded by
\begin{equation*}
C_1 ||(\hat{u}_0, \hat{\Theta}_0)||_{\gamma,\beta}\frac{(l+1)^{2/3} AD^l e^{-\beta|k|+\omega' p}}{(2l+1)(1+p^2)(1+|k|)^\gamma}\mathcal{Q}_{2l+2}(|\beta k|),
\end{equation*}
\begin{equation*}
\left|\frac{k_j}{l}\left(P_k(\hat{u}_{1,j}\hat{*}\hat{H}^{[l-1]}), \hat{u}_{1,j}\hat{*}\hat{S}^{[l-1]}\right)\right|\textnormal{ and }
\left|\frac{k_j}{l}\left(P_k(\hat{H}^{[l-1]}_j\hat{*}\hat{u}_{1}), \hat{H}^{[l-1]}_j\hat{*}\hat{\Theta}_1\right)\right|
\end{equation*}
are bounded by 
\begin{equation*}
C_1 ||(\hat{u}_1, \hat{\Theta}_1)||_{\gamma,\beta}\frac{l^{2/3} AD^{l-1}e^{-\beta|k|+\omega' p}}{l(2l-1)(1+p^2)(1+|k|)^\gamma}\mathcal{Q}_{2l}(|\beta k|),
\end{equation*}
and 
\begin{equation*}
|P_k(a e_2\hat{S}^{[l]})|\leq a\frac{e^{\omega' p}e^{-\beta |k|}AD^l}{(1+p^2)(1+|k|)^{\gamma}}\frac{\mathcal{Q}_{2l}(|\beta k|)}{(2l+1)^2}.
\end{equation*}
\noindent In the above, $C_1=C_1(d)$ is defined in \ref{6.10.}.
\end{Lemma}

\begin{pf} For the first inequality, we use (\ref{1-4.60}) and then apply \ref{6.10.} to get
\begin{align*}
(1+p^2)e^{-\omega' p}&|k_j\hat{u}_{0,j}\hat{*}(\hat{H}^{[l]}, \hat{S}^{[l]})|\\
&\leq||\hat{u}_0||_{\gamma,\beta}\frac{A D^l}{(2l+1)^2}|k|\int_{k'\in\mathbb{R}^d}\frac{e^{-\beta(|k'|+|k-k'|)}}{(1+|k'|)^{\gamma}(1+|k-k'|)^{\gamma}}\mathcal{Q}_{2l}(\beta|k'|)dk'\\ 
&\leq C_1(l+1)^{2/3}||\hat{u}_0||_{\gamma,\beta}\frac{A D^l e^{-\beta |k|}}{(2l+1)(1+|k|)^{\gamma}}\mathcal{Q}_{2l+2}(\beta|k|).
\end{align*}
\noindent The other inequalities are proved similarly and the last is simply the statement of the assumed bound.
\end{pf}

\begin{Lemma}\label{4.8.} Suppose $(\hat{H}^{[l]}, \hat{S}^{[l]})$ satisfies $(\ref{1-4.60})$ for $l\geq1$. Then 
\begin{equation*}
\left|\frac{k_j}{l}\left(P_k(\hat{H}_j^{[l-1]}(\cdot,0)\hat{*}\hat{H}^{[0]}(\cdot, p)),\hat{H}_j^{[l-1]}(\cdot,0)\hat{*}\hat{S}^{[0]}(\cdot, p)\right)\right|
\end{equation*}
is bounded by
\begin{equation*}
C_1\frac{(l+1)^{2/3}\tilde A^2 \tilde D^{l-1}e^{-\beta |k|+\alpha' p}}{l(2l-1)(1+|k|)^{\gamma}(1+p^2)}\mathcal{Q}_{2l}(\beta |k|).
\end{equation*}
\end{Lemma}
 \begin{pf} Using $(\ref{1-4.60})$ with $p=0$ and $(\ref{4.59})$ with $ A=3||(\hat{u}_1, \hat{\Theta}_1)||_{\gamma, \beta}$ along with \ref{6.10.}, we get
\begin{align*}
(1+p^2)e^{-\omega' p}&\left|\right.\frac{k_j}{l}[\hat{H}_j^{[l-1]}(\cdot,0)\hat{*}(\hat{H}^{[0]}, \hat{S}^{[0]})(\cdot, p)]\left.\right|\\ 
&\leq \frac{ A^2 D^{l-1}}{l(2l-1)^2}|k|\int_{k'\in \mathbb{R}^d}\frac{e^{-\beta(|k'|+|k-k'|)}}{(1+|k'|)^{\gamma}(1+|k-k'|)^{\gamma}}\mathcal{Q}_{2l-2}(\beta |k'|)dk'\\ 
&\leq C_1\frac{l^{2/3} A^2 D^{l-1}e^{-\beta |k|}}{l(2l-1)(1+|k|)^{\gamma}}\mathcal{Q}_{2l}(\beta |k|)
\end{align*}
\noindent From this the lemma follows after using Lemma \ref{2.3.}.
\end{pf}

\begin{Lemma}\label{4.9.} Suppose $(\hat{H}^{[l_1]}, \hat{S}^{[l_1]})$ and $(\hat{H}^{[l-l_1-1]}, \hat{S}^{[l-l_1-1]})$ satisfies $(\ref{1-4.60})$ for $l_1=1, \dots , l-2$ where $l\geq 2$. Then for $C_8=82$ and  $C_7=C_7(d)$ given in \ref{6.9.}, we have
\begin{equation*}
\left|k_j\sum_{l_1=1}^{l-2}\frac{l_1!(l-l_1-1)!}{l!}\left(P_k(\hat{H}^{[l_1]}_j(\cdot,0)\hat{*}\hat{H}^{[l-l_1-1]}(\cdot,p)), \hat{H}^{[l_1]}_j(\cdot,0)\hat{*}\hat{S}^{[l-l_1-1]}(\cdot,p)\right)\right|
\end{equation*}
is bounded by 
\begin{equation*}
C_8C_72^{\gamma}\pi A^2D^{l-1}\frac{e^{-\beta|k|+\omega' p}}{3\beta^d(1+p^2)(1+|k|)^{\gamma}}\frac{l\mathcal{Q}_{2l}(\beta |k|)}{(2l+3)^2}.
\end{equation*}
\end{Lemma}
\noindent The proof is the same as in \cite{smalltime} the only difference is a change in the constants arising when \ref{6.9.} in $\mathbb{R}^2$ or $\mathbb{R}^3$ is applied.

\begin{Lemma}\label{4.10.} Suppose $(\hat{H}^{[l]}, \hat{S}^{[l]})$ satisfies $(\ref{1-4.60})$ for $l\geq0$. Then
\begin{multline}\nonumber
\left|k_j\int_0^p\left(P_k(\hat{H}^{[l]}_j(\cdot,p-s)\hat{*}\hat{H}^{[0]}(\cdot,s)), \hat{H}^{[l]}_j(\cdot,p-s)\hat{*}\hat{S}^{[0]}(\cdot,s)\right)ds\right| \hspace{.5 in}\\ \nonumber
\leq C_1M_0A^2 D^{l}\frac{(l+1)^{2/3}e^{-\beta |k|+\omega' p}}{(2l+1)(1+|k|)^{\gamma}(1+p^2)}\mathcal{Q}_{2l+2}(\beta |k|).
\end{multline}
\noindent In the above, $M_0$, defined in Lemma \ref{2.6.}, is such that 
\begin{equation*}
\int_0^p\frac{1}{(1+(p-s)^2)(1+s^2)}ds\leq\frac{M_0}{1+p^2}.
\end{equation*}
\end{Lemma}

\begin{pf} Using (\ref{1-4.60}) for the first inequality and \ref{6.10.} and Lemma \ref{2.6.} for the second, we have
\begin{align}\nonumber
&\left|k_j\int_0^p\left(P_k(\hat{H}^{[l]}_j(\cdot,p-s)\hat{*}\hat{S}^{[0]}(\cdot,s)), (\hat{H}^{[l]}_j(\cdot,p-s)\hat{*}\hat{S}^{[0]}(\cdot,s))\right)ds\right|\leq \\ \nonumber
&|k|\frac{ A^2 D^l}{(2l+1)^2}\int_0^p\int_{k'\in\mathbb{R}^d} \frac{e^{-\beta|k'|+|k-k'|}e^{\omega'(p-s)+\omega's}}{(1+(p-s)^2)(1+s^2)(1+|k'|)^{\gamma}(1+|k-k'|)^{\gamma}}\mathcal{Q}_{2l}(\beta|k'|)ds dk'\\ \nonumber
&\qquad\qquad\leq C_1M_0 A^2 D^{l}\frac{(l+1)^{2/3}e^{-\beta |k|+\alpha p}}{(2l+1)(1+|k|)^{\gamma}(1+p^2)}\mathcal{Q}_{2l+2}(\beta |k|).
\end{align}
\end{pf}

\begin{Lemma}\label{4.11.} We have 
\begin{align}\nonumber \left|k_j\left(P_k(\hat{u}_{0,j}\hat{*}\hat{u}_1),\hat{u}_{0,j}\hat{*}\hat{\Theta}_1\right)\right.&+k_j\left.\left(P_k(\hat{u}_{1,j}\hat{*}\hat{u}_0), \hat{u}_{1,j}\hat{*}\hat{\Theta}_0\right)\right|\\ \nonumber
& \leq\frac{2C_0|k|e^{-\beta |k|}}{(1+|k|)^{\gamma}}||(\hat{u}_0,\hat{\Theta}_0)||_{\gamma, \beta}||(\hat{u}_1,\hat{\Theta}_1)||_{\gamma, \beta}\\ \nonumber
 \left|k_j\left(P_k(\hat{u}_{1,j}\hat{*}\hat{u}_1),\hat{u}_{1,j}\hat{*}\hat{\Theta}_1\right)\right|&\leq \frac{|k|e^{-\beta |k|}C_0}{(1+|k|)^{\gamma}}||\hat{u}_1, \hat{\Theta}_1||_{\gamma, \beta}^2.
\end{align}
\end{Lemma}

\begin{pf} This follows directly from Corollary \ref{2.2.} and Lemma \ref{2.3.}. 
\end{pf}

\begin{Lemma}\label{4.12.} For the case $l=1$, we have 
\begin{equation*}
|(\hat{H}^{[1]}, \hat{S}^{[1]})(k,p)|\leq \frac{e^{\omega' p}e^{-\beta |k|}AD}{(1+p^2)(1+|k|)^{\gamma}}\mathcal{Q}_{2}(|\beta k|),
\end{equation*}
\noindent where
\begin{align*} 
AD\geq &C_6\left(\frac{C_0}{\beta}\right. || (\hat{u}_0, \hat{\Theta}_0)||_{\gamma, \beta}||(\hat{v}_1, \hat{\Theta}_1)||_{\gamma, \beta}+M_1\frac{2}{\beta^2}||(\hat{u}_1, \hat{\Theta}_1)||_{\gamma, \beta}\\ 
 &\left.\qquad+C_1M_0A^2+2C_1A||(\hat{u}_0, \hat{\Theta}_0)||_{\gamma, \beta}+\frac{a A}{4}\right).
\end{align*}
\end{Lemma}
\begin{pf} Lemma \ref{4.4.} with $l=0$ tells us that 
\begin{equation*}
|(\hat{H}^{[1]}, \hat{S}^{[1]})(k,p)|\leq C_6 \sup_{p'\in [0,p]}|(\hat{R}^{[0]}_1, \hat{R}^{[0]}_2)(k,p')|
\end{equation*}
\noindent since $(\hat{H}^{[0]}, \hat{S}^{[0]})(k,0)=0$. We use Lemma \ref{4.6.}, Lemma \ref{4.10.}, and Lemma \ref{4.11.} to bound the terms appearing in $R_m$s. 
\begin{align*}
&|(\hat{R}^{[0]}_1, \hat{R}^{[0]}_2)(k,p)|\leq \frac{2C_0|k|e^{-\beta |k|}}{(1+|k|)^{\gamma}}||(\hat{u}_0, \hat{\Theta}_0)||_{\gamma, \beta}||(\hat{u}_1, \hat{\Theta}_1)||_{\gamma, \beta}\\ 
&\qquad +M_1\frac{|k|^2e^{-\beta |k|}}{(1+|k|)^{\gamma}}||(\hat{v}_1, \hat{\Theta}_1)||_{\gamma,\beta}+C_1M_0A^2 \frac{e^{-\beta |k|
+\omega' p}}{(1+|k|)^{\gamma}(1+p^2)}\mathcal{Q}_{2}(\beta |k|)\\ 
&\qquad +2C_1||(\hat{u}_0, \hat{\Theta}_0)||_{\gamma,\beta}\frac{A e^{-\beta|k|+\omega' p}}{(1+p^2)(1+|k|)^\gamma}\mathcal{Q}_{2}(|\beta k|)+a\frac{e^{\omega' p}e^{-\beta |k|}A}{(1+p^2)(1+|k|)^{\gamma}}
\end{align*}
\noindent The lemma now follows since $4|k|\leq \frac{2\mathcal{Q}_2}{\beta}$ and $|k|^2\leq\frac{2\mathcal{Q}_2}{\beta^2}$.
\end{pf}

{\bf Proof of Lemma \ref{4.2.}} Lemma \ref{4.12.} and $(\ref{4.59})$ prove the base case. Suppose, for the purpose of induction, that for $l\geq1$ $(\ref{1-4.60})$ holds. Then by Lemma \ref{4.4.} we need only prove a bound for $|(\hat{R}_1^{[l]}, \hat{R}_2^{[l]})|$ whose terms we bounded in the previous lemmas. 
\begin{align*}
&|(\hat{R}_1^{[l]}, \hat{R}_2^{[l]})|\leq \frac{AD^{l-1}e^{-\beta |k|+\omega' p}}{(2l+3)^2(1+p^2)(1+|k|)^{\gamma}}\mathcal{Q}_{2l+2}(\beta |k|)\left\{\frac{C_1M_0AD(l+1)^{2/3}(2l+3)^2}{(2l+1)}\right. \\ \nonumber
&+\frac{C_1A(l+1)^{2/3}(2l+3)^2}{4l(2l-1)}+\frac{C_8C_72^{\gamma}\pi Al}{12\beta^d}+\frac{C_1l^{2/3}||(\hat{u}_1, \hat{\Theta}_1)||_{\gamma,\beta}(2l+3)^2}{2l(2l-1)}\\ 
&\left. +2C_1D||(\hat{u}_0, \hat{\Theta}_0)||_{\gamma,\beta}\frac{(l+1)^{2/3}(2l+3)^2}{2l+1} +25\delta_{l,1}\frac{C_0}{A\beta}||(\hat{u}_1, \hat{\Theta}_1)||^2_{\gamma,\beta}+\frac{a D(2l+3)^2}{4(2l+1)^2}\right\}.
\end{align*} 
\noindent We also note that as $(\hat{H}^{[l]}, \hat{S}^{[l]})$ satisfies $(\ref{1-4.60})$,
\begin{align*}
\frac{|k|^2|(\hat{H}^{[l]}, \hat{S}^{[l]})(k,0)|}{(l+1)(l+2)}&\leq \frac{|k|^2e^{-\beta |k|}AD^l \mathcal{Q}_{2l}(\beta|k|)}{(l+1)(l+2)(1+|k|)^{\gamma}(2l+1)^2}\\ 
&\leq\frac{AD^{l}e^{-\beta |k|+\alpha' p}}{(2l+3)^2(1+p^2)(1+|k|)^{\gamma}}\mathcal{Q}_{2l+2}(\beta |k|)\frac{6}{\beta^2}.
\end{align*}
\noindent Here, we used the following two facts
\begin{equation*}
\frac{y^2\mathcal{Q}_{2l}(y)}{(2l+2)(2l+1)}\leq \mathcal{Q}_{2l+2}(y) \textnormal{ and }\frac{(2l+2)(2l+3)^2}{(l+1)(l+2)(2l+1)}\leq 6.
\end{equation*}
\noindent Thus, for $D$ chosen, independently of $l,k,$ and $p$, large enough so
\begin{align*}
D^2\geq& C_6\left\{\frac{C_1M_0AD(2l+3)^2}{(l+1)(2l+1)}+\frac{C_1A(2l+3)^2}{4(l+1)l(2l-1)}+\frac{C_8C_72^{\gamma}\pi Al}{12\beta^d(l+1)^{5/3}}\right. \\ 
& \qquad+\frac{C_1||(\hat{u}_1, \hat{\Theta}_1)||_{\gamma,\beta}(2l+3)^2}{2(l+1)^{5/3}l^{1/3}(2l-1)} +2C_1D||(\hat{u}_0, \hat{\Theta}_0)||_{\gamma,\beta}\frac{((2l+3)^2}{(l+1)(2l+1)}\\ 
&\qquad +25\delta_{l,1}\frac{C_0}{A2^{5/3}\beta}||(\hat{u}_1, \hat{\Theta}_1)||^2_{\gamma,\beta}+\left.\frac{a D(2l+3)^2}{4(l+1)^{5/3}(2l+1)}\right\}+M_1\frac{6D}{\beta^2},
\end{align*}
\noindent $(\ref{1-4.60})$ holds and the lemma is proved. 

As $\mathcal{Q}_{2l}(\beta|k|)\leq 4^l e^{|\beta k|/2}$,
\begin{equation}\label{series}
(\hat{H}, \hat{S})(k,p;p_0)=\sum_{l=0}^{\infty}(\hat{H}^{[l]}, \hat{S}^{[l]})(k,p_0)(p-p_0)^l
\end{equation}
\noindent is convergent for $|p-p_0|\leq\frac{1}{4D}$ where $D$ is independent of $p_0$. Moreover, the following lemma proved in \cite{smalltime} says that these series are indeed 
local representations of the solution $(\hat{H}, \hat{S})(k,p)$.

\begin{Lemma}\label{4.13.} The unique solution to ($\ref{1-2.9}$) satisfying $(\hat{H}, \hat{S})(k,0)=0$ guaranteed in Lemma \ref{2.8.} has a local representation given by $(\hat{H}, \hat{S})(k,p;p_0)$ for $p_0\in \mathbb{R}^+$. So, the solution is analytic on $\mathbb{R}^{+}\cup\{ 0\}$. 
\end{Lemma}

{\bf Proof of Theorem \ref{Borel summability} i)} Using Lemma \ref{4.2.} and the fact that $||g||_{L^{\infty}}\leq ||\hat{g}||_{L^1}$ we know that
\begin{align*}
|(H^{[l]},S^{[l]})(x,p_0)|&\leq\frac{8\pi A(4B)^le^{\omega p_0}}{\beta(2l+1)^2(1+p_0^2)}\\
|D(H^{[l]},S^{[l]})(x,p_0)|&\leq\frac{8\pi A(4B)^le^{\omega p_0}}{\beta(2l+1)^2(1+p_0^2)}\\
|D^2(H^{[l]},S^{[l]})(x,p_0)|&\leq\frac{16\pi A(4B)^le^{\omega p_0}}{\beta^2(2l+1)^2(1+p_0^2)}
\end{align*}
and the series (\ref{series}) converges for $|p-p_0|<\frac{1}{4B}$. By Lemma \ref{4.13.} the series is the 
local representation of the solution guaranteed to exist by Lemma \ref{2.8.} which is zero at $p=0$. Combining this with the facts that the solution is analytic in a neighborhood of zero and exponentially bounded for large p, recall ($\hat{H},\hat{S}\in \mathcal{A}^{\omega}$), implies Borel summability in $1/t$. Watson's Lemma then implies as $t\rightarrow 0^+$
\begin{equation*} 
(u, \Theta)(x,t)\sim (u_0, \Theta_0)(x)+
\sum_{m=1}^{\infty}(u_m, \Theta_m) (x)t^m
\end{equation*}
\noindent where $|(u_m, \Theta_m)(x)|\leq m! A_0 D_0^m$ with constants
$A_0$ and $D_0$ generally dependent on the initial condition and forcing through Lemma \ref{3.2.}.
\section{Extension of Existence Time}\label{sec7}

\indent We have shown by Theorem \ref{existence} that there is a unique solution to (\ref{1-2.18}) within the class of locally integrable functions, which are exponentially bounded in $p$, uniformly in $x$. Further, the solution $(\hat{H}, \hat{S})(k,p)$ generates a smooth solution to the Boussinesq equation for $t\in [0,\omega^{-1})$ where $\omega$ is the exponential growth rate of the integral equation (\ref{1-2.18}), and we showed that the solution is Borel summable. The question of global existence is then reduced to a question of exponential growth for the integral equation solution. If $(\hat{H}, \hat{S})(k,p)$ grows subexponentially, then global existence follows. The exponential growth rate $\omega$ previously found is suboptimal and ignores possible cancellations in the integrals. If we improve the estimates, we get a longer interval of existence. Here we present two examples of cases which can result in longer interval of existence.

\subsection{Improved Radius of Convergence}\label{subsec1}

When the initial data and forcing are analytic Borel summability given in Theorem \ref{Borel summability} implies that 
\begin{equation}\label{4.20}
(\hat{H}, \hat{S})(k,p)=\sum_{m=1}^{\infty}(\hat{u}^{[m]}, \hat{\Theta}^{[m]})(k)\frac{p^{m-1}}{(m-1)!}=\sum_{m=0}^{\infty}(\hat{u}^{[m+1]}, \hat{\Theta}^{[m+1]})(k)\frac{p^{m}}{m!}
\end{equation}
\noindent has a finite radius of convergence depending on the size of the initial data and forcing. However, in the special case when the initial data and forcing have only a finite number of Fourier modes the radius of convergence is in fact independent of the size of the initial data or $f$. The argument allows forcing to be time dependent. 

\noindent{\bf Proof of Theorem \ref{Borel summability} ii)} For small time 
\begin{align*}
(u, \Theta)(k,t)&=(\hat{u}^{[0]}, \hat{\Theta}^{[0]})(k)+\sum_{m=1}^{\infty}(\hat u^{[m]}, \hat{\Theta}^{[m]})(k)t^m\\
\hat{f}(k,t)&=\hat{f}^{[0]}+\sum_{m=1}^{\infty}\hat{f}^{[m]}(k)t^m,
\end{align*}
\noindent where by $(\ref{FB2})$ for $m\geq 0$
\begin{align}\label{4.21}
(m+1)\hat{u}^{[m+1]}&=\hat{f}^{[m]}-\nu|k|^2\hat{u}^{[m]}-i k_j P_k\left(\sum_{l=0}^m\hat{u}_j^{[l]}\hat{*}\hat{u}^{[m-l]}\right)+a P_k (e_2\hat{\Theta}^{[m]})\\ \nonumber
(m+1)\hat{\Theta}^{[m+1]}&=-\mu|k|^2\hat{\Theta}^{[m]}-i k_j\left(\sum_{l=0}^m\hat{u}_j^{[l]}\hat{*}\hat{\Theta}^{[m-l]}\right).
\end{align}
Suppose the initial data and forcing have a finite number of Fourier modes. Let $K_1=\max(\sup_{k\in supp(\hat{u}^{[0]}, \hat{\Theta}^{[0]})}|k|,\, \sup_{k\in supp(\hat{f})}|k|)$. Then by induction on $k$ we have $\sup_{k\in supp(\hat{u}^{[m]}, \hat{\Theta}^{[m]})}|k|\leq (m+1)K_1$. Taking the $||\cdot||_{\gamma,\beta}$ norm of both sides of (\ref{4.21}) with respect to $k$ and writing
\begin{equation*}
a_m=||(\hat{u}^{[m]}, \hat{\Theta}^{[m]})||_{\gamma,\beta},\qquad b_m=||\hat{f}^{[m]}||_{\gamma,\beta},
\end{equation*}
we obtain
\begin{align*}
a_{m+1}\leq&\frac{1}{m+1}\left[b_m+M_1\left\||k|^2|(\hat{u}^{[m]},\hat{\Theta}^{[m]})|\right\|_{\gamma,\beta}\right.\\ \nonumber
&\left. \qquad \qquad +\sum_{l=0}^m\left\| |k||\hat{u}^{[l]}|\hat{*}|(\hat{u}^{[m-l]},\hat{\Theta}^{[m-l]})|\right\|_{\gamma,\beta}+a a_m\right]\\ 
\leq& \frac{b_m}{m+1}+\frac{a a_m}{m+1}+K_1^2M_1(m+1)a_m+2K_1C_0\sum_{l=0}^m a_la_{m-l}.
\end{align*}
\noindent Consider the formal power series $y_0(t):=\sum_{m=1}^{\infty}\tilde{a}_m t^m$, where $\tilde{a}_0=a_0$ and
\begin{equation}\label {4.23}
\tilde{a}_{m+1}=\frac{b_m}{m+1}+\frac{a\tilde{a}_m}{m+1}+K_1^2M_1(m+1)\tilde{a}_m+2K_1C_0\sum_{l=0}^m \tilde{a}_l\tilde{a}_{m-l}.
\end{equation}
\noindent Clearly, $a_m\leq\tilde{a}_m$, so $y_0(t)$ majorizes $||(\hat{u}, \hat{\Theta})(\cdot,t)||_{\gamma,\beta}$. If we multiply both sides of (\ref{4.23}) by $t^m$ and sum over $m$, then
\begin{equation*}
\sum_{m=0}^{\infty}\tilde{a}_{m+1}t^m=\sum_{m=0}^{\infty}\left(\frac{b_m+a \tilde{a}_m}{m+1} +K_1^2M_1(m+1)\tilde{a}_m+2K_1C_0\sum_{l=0}^m \tilde{a}_l\tilde{a}_{m-l}\right)t^m.
\end{equation*}
\noindent In other words, $y_0(t)$ is a formal power series solution to 
\begin{equation*}
\frac{1}{t}(y-\tilde{a}_0)=w+\frac{a}{t}\int_0^t y(\tau)d\tau+K_1^2M_1(ty)'+2K_1C_0 y^2, 
\end{equation*}
\noindent where $w(t)=\sum_{m=0}^{\infty}\frac{b_m}{m+1}t^m $. With the change of variables $s=1/t$, we have
\begin{equation*}
-K_1^2M_1y'+2K_1C_0s^{-1}y^2+(K_1^2M_1s^{-1}-1)y+(s^{-1}w+\tilde{a}_0)+a s\int_0^{1/s}y(\tau)d\tau=0.
\end{equation*}
\noindent A singularity of $B(y(s))$ in the Borel plane exhibits itself as an exponential small correction to $y_0$. So, we let $y=y_0+\delta$ and construct the equation for $\delta$:
\begin{equation*}
-K_1^2M_1\delta'+2K_1C_0s^{-1}(\delta^2+2y_0\delta)+(K_1^2M_1s^{-1}-1)\delta+a s\int_0^{1/s}\delta(\tau)d\tau=0.
\end{equation*}
\noindent If we assume $\delta$ is exponentially small, then to leading order the equation is
\begin{equation*}
-K_1^2M_1\delta'+\left[(4K_1C_0s^{-1}\tilde{a}_0+(K_1^2M_1)s^{-1}-1\right]\delta=0,
\end{equation*}
\noindent which yields
\begin{equation*}
\delta\sim e^{-K_1^{-2}M_1^{-1}s}s^{4\tilde{a}_0C_0 K_1^{-1}M_1^{-1}+1}.
\end{equation*}
So, the radius of convergence of $B(y)$ is at least $K_1^{-2}M_1^{-1}$ which is independent of the size of initial data as claimed. As $y$ majorizes our solution $(\hat{u}, \hat{\Theta})(k,t)$ the radius of convergence of $(\ref{4.20})$ is independent of the size of initial data or forcing as well.

\subsection{Improved Growth Estimates Based on Knowledge of the Solution in [0, p0].}

Let $(\hat{H},\hat{S})(k,p)$ be the solution to (\ref{1-2.18}) provided by Theorem \ref{existence}. Recall the definitions of $(\hat{H},\hat{S})^{(a)}$ and $(\hat{H},\hat{S})^{(s)}$ given by (\ref{beginbehave}) and (\ref{smallbehave}) and the functionals in (\ref{13.10}) and (\ref{13.11}). 
\noindent Now, let $(\hat{H}, \hat{S})^{(b)}=(\hat{H}, \hat{S})-(\hat{H}, \hat{S})^{(a)}$. It is convenient to write the integral equation for $(\hat{H}, \hat{S})^{(b)}$ for $p>p_0$,
\begin{align}\label{l8.63}
\hat{H}^{(b)}(k,p)&=\frac{\pi}{z}\int_{p_0}^p \mathcal{G}(z,z')\left(i k_j\right.\hat{G}_j^{[1],(b)}(k,p')+P_k[e_2\hat{s}^{(b)}(k,p'\left.)]\right)dp'+\hat{H}^{(s)}(k,p)\\ \nonumber
\hat{S}^{(b)}(k,p)&=\frac{i k_j\pi}{2|k|\sqrt{\mu p}}\int_{p_0}^p \mathcal{G}(\zeta ,\zeta ')\hat{G}_j^{[2],(b)}(k,p')dp'+\hat{S}^{(s)}(k,p),
\end{align}
\noindent where 
\begin{align*}
\hat{G}_j^{[1],(b)}(k,p)=-P_k[\hat{u}_{0,j}\hat{*}\hat{H}^{(b)}+\hat{H}_j^{(b)}\hat{*}\hat{u}_0+\hat{H}_j^{(a)}\, ^{\ast}_{\ast}\hat{H}^{(b)}+\hat{H}_j^{(b)}\, ^{\ast}_{\ast}\hat{H}^{(a)}+\hat{H}_j^{(b)}\, ^{\ast}_{\ast}\hat{H}^{(b)}]\\ 
\hat{G}_j^{[2],(b)}(k,p)=-[\hat{u}_{0,j}\hat{*}\hat{S}^{(b)}+\hat{H}_j^{(b)}\hat{*}\hat{\Theta}_0+\hat{H}_j^{(a)}\, ^{\ast}_{\ast}\hat{S}^{(b)}+\hat{H}_j^{(b)}\, ^{\ast}_{\ast}\hat{S}^{(a)}+\hat{H}_j^{(b)}\, ^{\ast}_{\ast}\hat{S}^{(b)}].
\end{align*}
\noindent We also define 
\begin{equation}
\hat{R}^{(b)}(k,p)=i k_j(\hat{G}_j^{[1]}, \hat{G}_j^{[2]})^{(b)}(k,p)+a P_k[e_2\hat{S}^{(b)}(k,p)].
\end{equation}

\noindent{\bf Proof of Theorem \ref{improved existence}}
We note that
\begin{multline*}
|R^{(b)}(k,p)|\leq \left(|k|\left[|\hat{u}_0|\hat{*}|(\hat{H}, \hat{S})^{(b)}|+|\hat{H}^{(b)}|\hat{*}|(\hat{u}_0, \hat{\Theta}_0)|+2|(\hat{H}, \hat{S})^{(a)}|\, ^{\ast}_{\ast}|(\hat{H}, \hat{S})^{(b)}|\right.\right.\\ 
\left.\left.+|\hat{H}^{(b)}|\, ^{\ast}_{\ast}|(\hat{H}, \hat{S})^{(b)}|\right]+a|\hat{H}^{(b)}|\right)(k,p),
\end{multline*}
\noindent where $|\cdot|$ is the usual euclidean norm. Let $\psi(p)=||(\hat{H}, \hat{S})^{(b)}(\cdot,p)||_{\gamma, \beta}$. Then
\begin{multline*}
\left\|\left(\frac{\mathcal{G}(z,z')}{z}(i k_j(\hat{G}_j^{[1]})^{(b)}(k,p)+a P_k[e_2\hat{S}^{(b)}(k,p)]), \frac{\mathcal{G}(\zeta,\zeta')}{\zeta}i k_j(\hat{G}_j^{[2]})^{(b)}(k,p)\right)\right\|_{\gamma,\beta} \\ 
\leq \mathcal{B}_0(k)\cdot\left(|k|\left[||\hat{u}_0||_{\gamma,
\beta}\psi(p)+\psi(p)||(\hat{u}_0, \hat{\Theta}_0)||_{\gamma,\beta}+ 2||(\hat{H}, \hat{S})^{(a)}||_{\gamma,\beta}*\psi(p)\right.\right.\\ 
\left.\left.+\psi(p)*\psi(p)\right]+a\psi(p)\right)(k,p)= \left(\mathcal{B}_1\psi+\mathcal{B}_2*\psi+\mathcal{B}_3\psi*\psi+\mathcal{B}_4\psi\right)(p).
\end{multline*}
\noindent Taking the $(\gamma,\beta)$ norm in $k$ on both sides of (\ref{l8.63}) and multiplying by $e^{-\omega p}$ for $\omega\geq\omega_0\geq0$ and integrating from $p_0$ to $M$ gives
\begin{multline*}
L_{p_0,M}:=\int_{p_0}^Me^{-\omega p}\psi(p)dp\leq \int_{p_0}^Me^{-\omega p}\int_{p_0}^p\left(\mathcal{B}_1\psi+ \mathcal{B}_2*\psi+\mathcal{B}_3\psi*\psi+\mathcal{B}_4\psi\right)(p')dp'dp\\  +\int_{p_0}^Me^{-\omega p}\psi^{(s)}(p)dp \leq \int_{p_0}^M\int_{p'}^{M}e^{-\omega (p-p')}e^{-\omega p'}\left(\mathcal{B}_1\psi+ \mathcal{B}_2*\psi+\mathcal{B}_3\psi*\psi+\mathcal{B}_4\psi\right)(p')dp dp'\\ 
+\int_{p_0}^Me^{-\omega p}\psi^{(s)}(p)dp\leq \frac{1}{\omega}\int_{p_0}^Me^{-\omega p'}\left(\mathcal{B}_1\psi+ \mathcal{B}_2*\psi+\mathcal{B}_3\psi*\psi+\mathcal{B}_4\psi\right)(p')dp'\\ 
+\int_{p_0}^Me^{-\omega p}\psi^{(s)}(p)dp,
\end{multline*}
\noindent where $\psi^{(s)}=||(\hat{H}, \hat{S})^{(s)}(\cdot,p)||_{\gamma,\beta}$. Recalling that $\psi=0$ on $[0,p_0]$, we note that for any $u$
\begin{equation*}
\int_{p_0}^Me^{-\omega p}(\psi *u)(p)dp=\int_{p_0}^M\psi(s)e^{-\omega s}\int_{0}^{M-s}e^{-\omega p}u(p)dp ds.
\end{equation*}
\noindent Using this, we obtain
\begin{multline}\nonumber
L_{p_0,M}\leq\frac{1}{\omega}\left\{(\mathcal{B}_1+\int_0^{M-p_0}e^{-\omega p}\mathcal{B}_2(p)dp)L_{p_0,M}+\mathcal{B}_3L_{p_0,M}^2+\mathcal{B}_4L_{p_0,M}\right\}+b\omega^{-1}\\ \nonumber
\leq \omega ^{-1}\left\{\epsilon_1L_{p_0,M}+\mathcal{B}_3L_{p_0,M}^2\right\}+b\omega^{-1}.
\end{multline}
\noindent For 
\begin{equation}\nonumber
\epsilon_1<\omega\quad \mbox{and}\quad(\epsilon_1-\omega)^2>4\mathcal{B}_3b,
\end{equation}
\noindent we get an estimate for $L_{p_0,M}$ that is independent of $M$. Namely,
\begin{equation*}
L_{p_0,M}\leq\frac{1}{2\mathcal{B}_3}\left[\omega-\epsilon_1-\sqrt{(\epsilon_1-\omega)^2-4\mathcal{B}_3b}\right].
\end{equation*}

So, $||(\hat{H}, \hat{S})(\cdot,p)||_{\gamma,\beta}\in L^1(e^{-\omega p}dp)$, and the solution to the Boussinesq exists for $t\in(0,\omega^{-1})$ for $\omega$ sufficiently large so that 
\begin{equation*}
\omega\geq \omega_0 \quad\mbox{and}\quad \omega>\epsilon_1+2\sqrt{\mathcal{B}_3b}.
\end{equation*}
Equivalently, we could choose our original $\omega_0$ large enough so that $ \omega_0>\epsilon_1+2\sqrt{\mathcal{B}_3b}$. This completes the proof of Theorem \ref{improved existence}.


\appendix
\section{Bessel Function Representation of the Kernel} 
\begin{Lemma}\label{kernel} The kernel $\mathcal{G}(z,z')$ given by
\begin{equation*}
\mathcal{G}(z,z')=z'(-J_1(z)Y_1(z')+Y_1(z)J_1(z')),\textnormal{where } z=2|k|\sqrt{\nu p} \textnormal{ and } z'=2|k|\sqrt{\nu p'}
\end{equation*}
satisfies $\frac{\pi}{z}\mathcal{G}(z,z')=\mathcal{H}^{(\nu)}(p,p',k)$ with $\mathcal{H}^{(\nu)}$ given by (\ref{Chap2H}).
\end{Lemma}
\textit{Proof.} We will show that $\mathcal{H}^{(\nu)}(p,p',k)$ solves $(p\partial_{pp}+2\partial_{p}+\nu|k|^2)\mathcal{H}^{(\nu)}=0$ for $0<p'<p$ with the condition that $\mathcal{H}^{(\nu)}(p,p',k)\rightarrow 0$ and $\mathcal{H}^{(\nu)}_p(p,p',k)\rightarrow \frac{1}{p}$ as $p'$ approaches p from below.

First, we notice that
\begin{equation*}
\mathcal{H}^{(\nu)}(p,p',k)=\frac{p'}{p}\int_1^{p/p'}F(\eta)ds,
\end{equation*}
 where
\begin{equation*}
\eta=\nu |k|^2p\left(1-\frac{sp'}{p}\right)\left(1-\frac{1}{s}\right),\quad F(\eta)=\frac{1}{2\pi i}\int_{C}\zeta^{-1}e^{\zeta-\eta\zeta^{-1}}d\zeta,
\end{equation*}
and $C$ is the contour starting and $\infty e^{-\pi i}$ turning around the origin in counterclockwise direction and ending at $\infty e^{\pi i}$. In the appendix of \cite{longtime}, it is shown that $F$ is entire, $F(0)=1$, and $F$ satisfies $\eta F''(\eta)+F'(\eta)+F(\eta)=0$. We will use these facts as given. As $F$ is continuous and the interval of integration shrinks to length zero, $\mathcal{H}^{(\nu)}(p,p',k)\rightarrow 0$ as $p'$ tends to $p$ from below. For $p> p'$, $\mathcal{H}^{(\nu)}$ is twice differentiable in $p$ as $F$ is twice continuously differentiable. Moreover, we have
\begin{align*}
\mathcal{H}^{(\nu)}_p(p,p',k)&=-\frac{1}{p}\mathcal{H}^{(\nu)}(p,p',k)+\frac{1}{p}F(0)+\frac{p'}{p}\int_1^{p/p'}F'(\eta)\frac{d\eta}{dp}ds,\\ 
(p\mathcal{H}^{(\nu)}_p)_p&=-\mathcal{H}^{(\nu)}_p+F'(0)\nu|k|^2(1-\frac{p'}{p})+p'\int_1^{p/p'}F''(\eta)\left(\frac{d\eta}{dp}\right)^2ds,
\end{align*}
where the second equality uses that $\frac{d\eta}{dp}=\nu|k|^2\left(1-\frac{1}{s}\right)$ is $p$ independent. Thus, as $F(0)=1$, we have $\mathcal{H}^{(\nu)}_p(p,p',k)\rightarrow \frac{1}{p}$ as $p'$ tends to $p$ from below.
We notice that 
\begin{equation*}
\left(\frac{d\eta}{dp}\right)^2=\frac{\eta\nu|k|^2}{p}-\frac{\nu|k|^2(s-1)}{p}\frac{d\eta}{ds}.
\end{equation*}
So, integrating by parts and using $\eta F''(\eta)+F'(\eta)+F(\eta)=0$, we have
\begin{align*}
(p\mathcal{H}^{(\nu)}_p)_p&+\mathcal{H}^{(\nu)}_p=F'(0)\nu|k|^2(1-\frac{p'}{p})+p'\int_1^{p/p'}F''(\eta)\left(\frac{\eta\nu|k|^2}{p}\right)ds\\ 
&\qquad \qquad-p'\int_1^{p/p'}\frac{d}{ds}(F'(\eta))\frac{\nu|k|^2(s-1)}{p}ds\\ 
=&\frac{\nu|k|^2p'}{p}\int_1^{p/p'}\eta F''(\eta)ds+\frac{p'\nu|k|^2}{p}\int_1^{p/p'}F'(\eta)ds=-\nu|k|^2\mathcal{H}^{(\nu)}.
\end{align*}
In other words, $p\mathcal{H}^{(\nu)}_{pp}+2\mathcal{H}^{(\nu)}_p+\nu|k|^2\mathcal{H}^{(\nu)}=0$, and the lemma is proved.

\begin{Lemma}\label{U_0} We also have the representation in terms of Bessel functions
\begin{equation*}
\mathcal{L}^{-1}\left(\frac{1-e^{-\nu|k|^2\tau^{-1}}}{\nu|k|^2}\right)(p)=\frac{2J_1(z)}{z}.
\end{equation*}
\end{Lemma}
\begin{pf}
Notice that by contour deformation the contribution from $\frac{1}{\nu|k|^2}$ is zero. Factoring out $|k|\sqrt{\nu p}$ in the exponent and using the change of variables $\frac{\tau\sqrt{p}}{|k|\sqrt{\nu}}\rightarrow w$, we have
\begin{equation*}
\mathcal{L}^{-1}\left(\frac{1-e^{-\nu|k|^2\tau^{-1}}}{\nu|k|^2}\right)(p)=\frac{-1}{2\pi i}\int_{c-i\infty}^{c+i\infty}\frac{e^{|k|\sqrt{\nu p}(w-w^{-1})}}{|k|\sqrt{\nu p}}dw
=2\frac{J_1(z)}{z}.
\end{equation*}
\end{pf}

\section{Fourier Inequalities in Two Dimensions}\label{sec8}

In the appendix of \cite{smalltime}, Fourier inequalities are developed in $\mathbb{R}^3$. We present the counterparts to those inequalities for $\mathbb{R}^2$ here. Where a lemma is referenced from this section, we use either the $\mathbb{R}^2$ or $\mathbb{R}^3$ version
as appropriate. The basic idea is that in 2-d \ref{6.6.} below differs by a constant from 3-d case. All other lemmas are basically the same for $\mathbb{R}^2$ or $\mathbb{R}^3$ once the change in \ref{6.6.} is taken into account.
\begin{Definition}\label{d6.1.} Define the polynomial
\begin{equation*}
P_n(z)=\sum_{j=0}^n\frac{n!}{j!}z^j.
\end{equation*}
\end{Definition}
\begin{Lemma}\label{6.5.} For all $y\geq 0$ and integers $n\geq m \geq 0$, we have 
\begin{equation*}
y^{m+1}\int_0^{\infty}e^{-y(\rho -1)[1+sgn(\rho -1)]}\rho ^m P_n(y|1-\rho|)d\rho\leq m!n!\mathcal{Q}_{m+n+1}(y).
\end{equation*}
\end{Lemma}
\noindent Proof can be found in \cite{smalltime}.

\begin{Proposition}\label{2-dbound} Let n be an integer no less than $0$ and $r\geq 0$ and $\rho \geq0$ fixed. Then
\begin{equation*}
\int_0^{2\pi}e^{-|\rho-re^{i \theta}|}|\rho-re^{i \theta}|^n d\theta\leq 2\pi e^3 e^{-|\rho-r|}P_n(|r-\rho|)+\frac{4e}{\rho}e^{-|\rho-r|}P_{n+1}(|\rho-r|).
\end{equation*}
\end{Proposition}
\begin{pf}
\noindent Case 1. Suppose $0\leq r\leq 2$. Then for all $\theta$, 
\begin{equation*}
|\rho-r|\leq|\rho-re^{i\theta}|\leq |\rho-r|+|r-re^{i\theta}|\leq|\rho-r|+4. 
\end{equation*}
We also notice, for $x\geq 0$,
\begin{equation}\label{case1}
(x+1)^n=\sum_{j=0}^nx^j\frac{n!}{j!(n-j)!}\leq\sum_{j=0}^nx^j\frac{n!}{j!}=P_n(x).
\end{equation}
Further, for $x,\, a\geq 0$,
\begin{equation}\label{P}
\sum_{j=0}^n\frac{(x+a)^j}{j!}=\sum_{m=0}^n\sum_{j=m}^n\frac{x^{j-m}a^m}{m!(j-m)!}\leq\sum_{m=0}^n\frac{a^m}{m!}\sum_{j=0}^{n}\frac{x^{j}}{m!}\leq e^a\sum_{j=0}^{n}\frac{x^{j}}{m!}.
\end{equation}
Thus, 
\begin{equation*}
|\rho-re^{i\theta}|^n\leq\sum_{j=0}^n(|\rho-r|+3)^j\frac{n!}{j!}\leq e^3P_n(|\rho-r|)
\end{equation*}
and
\begin{equation*}
\int_0^{2\pi}e^{-|\rho-re^{i \theta}|}|\rho-re^{i \theta}|^n d\theta\leq 2\pi e^3 e^{-|\rho-r|}P_n(|r-\rho|).
\end{equation*}
So, the proposition holds in this case.

\noindent Case 2. Suppose $r>2$. Let $\theta_1\in(0,\frac{\pi}{3})$ be such that $|r-re^{i\theta_1}|=1$. We split our integral into three pieces. For $\theta\in[0,\theta_1]$, 
\begin{equation*}
|\rho-r|\leq|\rho-re^{i\theta}|\leq |\rho-r|+|r-re^{i\theta}|\leq|\rho-r|+1. 
\end{equation*}
Applying (\ref{case1}) with $x=|\rho-r|$ gives, 
\begin{equation}\label{I1}
2\int_0^{\theta_1}e^{-|\rho-re^{i \theta}|}|\rho-re^{i \theta}|^n d\theta\leq\frac{2\pi}{3}e^{-|\rho-r|}P_n(|r-\rho|).
\end{equation}
Suppose $\theta\in[\theta_1,\pi-\theta_1]$. Let $z=|\rho-re^{i\theta}|=\sqrt{(\rho-r)^2+2\rho r(1-\cos\theta)}$. Then $d\theta=\frac{zdz}{\rho r \sin\theta}$. However, since $\theta\in[\theta_1,\pi-\theta_1]$,
\begin{equation}
\frac{1}{r\sin\theta}\leq\frac{1}{r\sin\theta_1}=\frac{1}{r\theta_1}\frac{\theta_1}{\sin\theta_1}.
\end{equation}
Now, notice that $\theta_1r\geq|r-re^{i\theta_1}|=1$ and $\frac{\theta_1}{\sin\theta_1}\leq\frac{\pi/3}{\sin(\pi/3)}< 2$ since $\theta_1\in[0,\frac{\pi}{3}]$. Hence,
\begin{equation*}
d\theta=\frac{zdz}{\rho r \sin\theta}\leq \frac{2zdz}{\rho}
\end{equation*} 
and
\begin{align}\label{I2}
2&\int_{\theta_1}^{\pi-\theta_1}e^{-|\rho-re^{i \theta}|}|\rho-re^{i \theta}|^n d\theta\leq\frac{4}{\rho}\int_{|\rho-re^{i\theta_1}|}^{|\rho-re^{i(\pi-\theta_1)}|}e^{-z}z^{n+1}dz\\ \nonumber
&=\frac{4}{\rho}\left(P_{n+1}(|\rho-re^{i\theta_1}|)e^{-|\rho-re^{i\theta_1}|}-P_{n+1}(|\rho-re^{i(\pi-\theta_1)}|)e^{-|\rho-re^{i(\pi-\theta_1)}|}\right).
\end{align}
We bound the positive contribution as in (\ref{P}) by
\begin{equation}\label{poscont}
\frac{4}{\rho}P_{n+1}(|\rho-re^{i\theta_1}|)e^{-|\rho-re^{i\theta_1}|}\leq\frac{4}{\rho}P_{n+1}(|\rho-r|+1)e^{-|\rho-r|}\leq\frac{4e}{\rho}P_{n+1}(|\rho-r|)e^{-|\rho-r|}.
\end{equation}
For $\theta\in[\pi-\theta_1, \pi]$, we again use (\ref{case1}) and get
\begin{equation*}
|\rho-re^{i\theta}|^n\leq(|\rho-re^{i(\pi-\theta_1)}|+1)^n\leq P_n(|\rho-re^{i(\pi-\theta_1)}|)
\end{equation*}  
and 
\begin{equation*}
|\rho-re^{i\theta}|=\sqrt{\rho^2-2\rho r\cos(\pi-\theta_1)+r^2+2\rho r(\cos(\pi-\theta_1)-\cos\theta)}\geq|\rho-re^{i(\pi-\theta_1)}|.
\end{equation*}
So,
\begin{equation}\label{I3}
2\int_{\pi-\theta_1}^{\pi}e^{-|\rho-re^{i \theta}|}|\rho-re^{i \theta}|^n d\theta\leq\frac{2\pi}{3}e^{-|\rho-re^{i(\pi-\theta_1)}|}P_n(|\rho-re^{i(\pi-\theta_1)}|).
\end{equation}
Now, we notice that $|\rho-re^{i(\pi-\theta_1)}|>\rho$, so
\begin{align*}
\frac{P_{n+1}(|\rho-re^{i(\pi-\theta_1)}|)}{\rho}&\geq(n+1)!\left(\frac{1}{\rho}+\sum_{j=1}^{n+1}\frac{|\rho-re^{i(\pi-\theta_1)}|^{j-1}}{j!}\right)\\
&\geq n!\sum_{j=0}^n\frac{n+1}{j+1}\frac{|\rho-re^{i(\pi-\theta_1})|^{j}}{j!}\geq P_n(|\rho-re^{i(\pi-\theta_1)}|).
\end{align*}
Thus,
\begin{equation}\label{negcont}
e^{-|\rho-re^{i(\pi-\theta_1)}|}\left(-\frac{2}{\rho}P_{n+1}(|\rho-re^{i(\pi-\theta_1)}|)+\frac{\pi}{3}P_n(|\rho-re^{i(\pi-\theta_1)}|)\right)<0.
\end{equation}
Adding the contributions from (\ref{I1}), (\ref{I2}), and (\ref{I3}) and using (\ref{poscont}) and (\ref{negcont}) gives
\begin{equation}
2\int_0^{\pi}e^{-|\rho-re^{i \theta}|}|\rho-re^{i \theta}|^n d\theta\leq\frac{2\pi}{3}e^{-|\rho-r|}P_n(|r-\rho|)+\frac{4e}{\rho}e^{-|\rho-r|}P_{n+1}(|\rho-r|).
\end{equation}
As all values of $r$ fall into one of these three cases, the proposition is proved.
\end{pf}

\begin{Lemma}\label{6.6.} If $m$ and $n$ are integers no less than $-1$, then
\begin{equation*}
|q|\int_{q'\in \mathbb{R}^d} e^{|q|-|q'|-|q-q'|}|q'|^m|q-q'|^n dq'\leq C_7(d) \pi (m+1)!(n+1)!\mathcal{Q}_{m+n+3}(|q|),
\end{equation*}
\noindent where $C_7(2)=6\pi e^3+4 e$ and $C_7(3)=2$.
\end{Lemma}

\begin{pf} We note that we may assume without loss of generality that $m\leq n$ since a change of variables $q'\rightarrow q-q'$ switches the roles of $m$ and $n$. Write $q=\rho e^{i \phi}$, $q'=re^{i\varphi}$ and $\theta=\varphi-\phi$. Let I be the integral on the left hand side. Then switching to polar coordinates gives
\begin{equation*}
I=\rho\int_0^{\infty}\int_0^{2\pi}e^{\rho-r-|\rho-re^{i\theta}|}r^m|\rho-re^{i\theta}|^n r dr d\theta.
\end{equation*}
For $n\geq 0$, using Proposition \ref{2-dbound} above gives,
\begin{equation*}
I\leq \rho\int_0^{\infty}e^{\rho-r}r^{m+1} e^{-|\rho-r|}(2\pi e^3P_n(|\rho-r|)+4 e \frac{P_{n+1}(|\rho -r|)}{\rho}dr.
\end{equation*}
\noindent Now, we let $\tilde{\rho}=\frac{r}{\rho}$. Then $d\tilde{\rho}=\frac{dr}{\rho}$ and $-|\rho-r|=-\rho(\tilde{\rho}-1)sgn(\tilde{\rho}-1)$, so
\begin{equation*}
I\leq \rho^{m+3}\int_0^{\infty}e^{-\rho(\tilde{\rho}-1) (1+sgn(\tilde{\rho}-1))}\tilde{\rho}^{m+1}(2\pi e^3P_n(\rho|\tilde{\rho}-1|)+4e \frac{P_{n+1}(\rho|\tilde{\rho} -1|)}{\rho}d\tilde{\rho}.
\end{equation*}
\noindent Applying \ref{6.5.} gives
\begin{align*}
I&\leq 2\pi e^3\rho(m+1)!n!\mathcal{Q}_{m+n+2}(\rho)+4 e (m+1)!(n+1)!Q_{m+n+3}\\
&\leq (6\pi e^3+4 e)(m+1)!(n+1)!\mathcal{Q}_{m+n+3}(\rho),
\end{align*}
\noindent where the last inequality follows as $m\leq n$, so 
\begin{align*}
\rho\sum_{j=0}^{m+n+2}\frac{2^{m+n+2-j}\rho^{j}}{j!}&\leq\sum_{j=1}^{m+n+3}\frac{2^{m+n+3-j}\rho^j}{(j-1)!}\\
&\leq \mathcal{Q}_{m+n+3}(\rho)(m+n+3)\leq3(n+1)\mathcal{Q}_{m+n+3}(\rho).
\end{align*} 
For $n=m=-1$, we use a slightly different approach. Assuming $q$ is not zero, we split the integral over two regions, a ball of radius $3|q|/2$ centered at zero and its compliment. For the compliment region we have $|q-q'|\geq|q|/2$, so
\begin{align*}
|q|\int_{|q'|\geq 3|q|/2}e^{|q|-|q'|-|q-q'|}&\frac{1}{|q'||q-q'|}dq'\\
&\leq 2e^{|q|/2}\int_0^{2\pi}\int_{3|q|/2}^{\infty}e^{-r} drd\theta=4\pi e^{-|q|}\leq4\pi.
\end{align*}
For the interior region we have 
\begin{equation*}
|q|\int_{|q'|\leq 3|q|/2}e^{|q|-|q'|-|q-q'|}\frac{1}{|q'||q-q'|}dq'\leq |q|\int_{|q'|\leq 3|q|/2}\frac{1}{|q'||q-q'|}dq'.
\end{equation*}
We now note that $\int_{|q'|\leq 3|q|/2}\frac{1}{|q'||q-q'|}dq'$ is bounded. Without trying to be precise we can bound the integral by $13\pi$ by spitting the region into two disks of radius $|q|/2$ centered at $0$ and $q$ and the compliment, call the compliment $D$. We have
\begin{equation*}
\int_{|q'|\leq|q|/2}\frac{1}{|q'||q-q'|}dq'\leq\frac{2}{|q|}\int_{|q'|\leq|q|/2}\frac{1}{|q'|}dq'\leq 2\pi.
\end{equation*}
Similarly,
\begin{equation*}
\int_{|q'-q|\leq|q|/2}\frac{1}{|q'||q-q'|}dq'\leq 2\pi.
\end{equation*}
Finally,
\begin{equation*}
\int_{D}\frac{1}{|q'||q-q'|}dq'\leq\frac{4}{|q|^2}\int_{D}dq'\leq \frac{4}{|q|^2}\int_{|q'|\leq 3|q|/2}dq'\leq9\pi.
\end{equation*}
Thus,
\begin{equation*}
|q|\int e^{|q|-|q'|-|q-q'|}\frac{1}{|q'||q-q'|}dq'\leq 13\pi|q|+4\pi\leq 13\pi(|q|+2)=13\pi Q_1(|q|)
\end{equation*}
for all nonzero $q$. Hence, the lemma is proved with $C_7(2)=6\pi e^3+4 e$.
\end{pf}

\begin{Lemma}\label{6.7.} For any $\gamma \geq 1$ and nonnegative integers $m$ and $n$, we have
\begin{align*}
|k|\int_{k'\in \mathbb{R}^d} &\frac{e^{-\beta(|k'|+|k-k'|)}}{(1+|k'|)^{\gamma}(1+|k-k'|)^{\gamma}}(\beta|k'|)^m(\beta|k-k'|)^n dk'\\ 
&\quad \leq \frac{C_7 \pi 2^{\gamma}e^{-\beta|k|}m!n!}{\beta^d(1+|k|)^{\gamma}}(m+n+2)\mathcal{Q}_{m+n+2}(\beta|k|).
\end{align*}
\end{Lemma}
\begin{pf} The proof is the same as the proof for 3-d given in \cite{smalltime} after using our new bound in \ref{6.6.}. 
\end{pf}

\begin{Lemma}\label{6.8.} For any $\gamma \geq 2$ and $n\in \mathbb{N}-0$, we have 
\begin{align*}
|k|\int_{k'\in \mathbb{R}^d} &\frac{e^{-\beta(|k'|+|k-k'|)}}{(1+|k'|)^{\gamma}(1+|k-k'|)^{\gamma}}|\beta(k-k')|^n dk'\\ 
&\leq \frac{C_7\pi 2^{\gamma}e^{-\beta|k|}}{\beta^{d-1}(1+|k|)^{\gamma}}\left\{(n-1)!\mathcal{Q}_{n+1}(\beta|k|)+\frac{3(n+1)!(\beta|k|)^{2/3}}{2\beta ^{2/3}}\sum_{j=0}^{n+1}\frac{(\beta|k|)^j}{j!}\right\}.
\end{align*}
\end{Lemma}

\begin{pf} We split the region into two integrals $\int_{|k'|\leq |k|/2}+\int_{|k'|\geq |k|/2}$. In the outer region, we have $(1+|k'|)^{-\gamma}\leq2^{\gamma}(1+|k|)^{-\gamma}$, and in the inner, we have $(1+|k-k'|)^{-\gamma}\leq 2^{\gamma}(1+|k|)^{-\gamma}$. We use this and $\gamma \geq 2$ for the first inequality and \ref{6.6.} for the second to get a bound for the outer region
\begin{align*}
|k|\int_{|k'|\geq |k|/2}& \frac{e^{-\beta(|k'|+|k-k'|)}}{(1+|k'|)^{\gamma}(1+|k-k'|)^{\gamma}}|\beta(k-k')|^n dk'\hspace{2.5 in}\\ 
&\leq \frac{2^{\gamma}e^{-\beta|k|}}{\beta^{d-1}(1+|k|)^{\gamma}}|q|\int_{q'\in \mathbb{R}^d}e^{|q|-|q'|-|q-q'|}|q-q'|^{n-2} dq'\\ 
&\leq\frac{C_7 \pi 2^{\gamma }e^{-\beta|k|}}{\beta^{d-1}(1+|k|)^{\gamma}}(n-1)!\mathcal{Q}_{n+1}(|q|).
\end{align*}
\noindent In the inner region, we also use $(1+|k'|)^{-\gamma}\leq (|k'|)^{-2+2/3}$, a change to polar coordinates as in the proof of \ref{6.6.}, and integration by parts to get
\begin{align*}
|k|&\int_{|k'|\leq|k|/2}\frac{e^{-\beta(|k'|+|k-k'|)}}{(1+|k'|)^{\gamma}(1+|k-k'|)^{\gamma}}|\beta(k-k')|^n dk'\\ 
&\leq \frac{2^{\gamma}e^{-\beta|k|}}{\beta^{d-1+2/3}(1+|k|)^{\gamma}}|q|\int_{|q'|\leq|q|/2}e^{|q|-|q'|-|q-q'|}|q'|^{-2+2/3}|q-q'|^{n}dq'\\ 
&=\frac{  2^{\gamma}e^{-\beta|k|}}{\beta^{d-1+2/3}(1+|k|)^{\gamma}}\rho\int_0^{\rho/2}\int_0^{2\pi}e^{\rho-r-|\rho-re^{i\theta}|}|\rho-re^{i\theta}|^n r^{-2+2/3}rd\theta dr\\ 
&\leq \frac{2^{\gamma}e^{-\beta|k|}}{\beta^{d-1+2/3}(1+|k|)^{\gamma}}\int_0^{\rho/2}r^{-1+2/3}(2\pi e^3\rho P_n(|\rho-r|)+4 e P_{n+1}(|\rho-r|)dr\\ 
&\leq \frac{2^{\gamma}e^{-\beta|k|}}{\beta^{d-1+2/3}(1+|k|)^{\gamma}}  \left(2\pi e^3n!\rho^{1+2/3}\sum_{j=0}^n\frac{\rho^{j}}{j!}\int_0^{1}\tilde{r}^{-1+2/3}(1-\tilde{r})^j d\tilde{r}\right.\\
&\qquad \left.4 e(n+1)!\rho^{2/3}\sum_{j=0}^{n+1}\frac{\rho^{j}}{j!}\int_0^{1}\tilde{r}^{-1+2/3}(1-\tilde{r})^j d\tilde{r}\right)\\ 
&\leq \frac{2^{\gamma}e^{-\beta|k|}}{\beta^{d-1+2/3}(1+|k|)^{\gamma}} \frac{C_73}{2}\pi (n+1)!\rho^{2/3}\sum_{j=0}^{n+1}\frac{\rho^{j}}{j!}.
\end{align*}
\end{pf}

The proof of the remaining lemmas is the same in 2-d as in 3-d after the change in bound given in \ref{6.8.} and can be found in \cite{smalltime}. Whenever Lemma 6.8. is invoked in \cite{smalltime} the 2-d proofs use \ref{6.8.}. 

\begin{Lemma}\label{6.9.} For any $\gamma \geq 1$ and nonnegative integers $l_1, l_2 \geq 0$, we have 
\begin{multline*}
|k|\int_{k'\in \mathbb{R}^d} \frac{e^{\beta(|k|-|k'|-|k-k'|)}}{(1+|k'|)^{\gamma}(1+|k-k'|)^{\gamma}}\mathcal{Q}_{2l_1}(\beta|k'|)\mathcal{Q}_{2l_2}(\beta|k-k'|)dk'\hspace{1.5 in}\\ 
\leq \frac{C_7 \pi 2^{\gamma}e^{-\beta|k|}}{3\beta^d(1+|k|)^{\gamma}}(2l_1+2l_2+1)(2l_1+2l_2+2)(2l_1+2l_2+3)\mathcal{Q}_{2l_1+2l_2+2}(\beta|k|).
\end{multline*}
\end{Lemma}

\begin{Lemma}\label{6.10.} If $\gamma\geq 2$ and $l\geq 0$, then
\begin{align*}
\frac{|k|}{(l+1)^{2/3}}\int_{k'\in\mathbb{R}^d}&\frac{e^{-\beta(|k'|+|k-k'|)}}{(1+|k'|)^{\gamma}(1+|k-k'|)^{\gamma}}\mathcal{Q}_{2l}(|\beta(k-k')|)dk'\\
&\leq \frac{C_1e^{-\beta|k|}}{(1+|k|)^{\gamma}}(2l+1)\mathcal{Q}_{2l+2}(\beta|k|),
\end{align*}
\noindent where
\begin{equation*}
C_1=C_1(d)=6C_7\pi2^{\gamma}\beta^{-d+1/3}+C_7\pi2^{\gamma}\beta^{-d+1}+\frac{1}{2}C_0\beta^{-1}.
\end{equation*}
\end{Lemma}





\begin{thebibliography}{999}
\bibitem{handbook} M. Abramowitz, I.A. Stegun, Handbook of mathematical functions with formulas, graphs, and mathematical tables, New York: Wiley-Interscience, 1970, pp 365 (see Formula 9.3.35-9.3.38).

\bibitem{Balser} W. Balser, Divergent Solutions of the heat equation: on an article of Lutz, Miyake
and Schafke, Pacific J. Math. 188 (1) (1999) 53-63. 

\bibitem{Braaksma} B.L.J. Braaksma, Transseries for a class of nonlinear Difference Equations,
J. Differ. Equations Appl. 7 (5) (2001) 717-750.

\bibitem{Cannon} J.R. Cannon, E. Dibenedetto, The Initial Value Problem for the Boussinesq Equation with Data in $L^p$, Lecture Notes in Mathematics: Approximation Methods for Navier-Stokes Problems, 771 (1979) 129-144. 

\bibitem{Costin-Costin} O. Costin, R.D. Costin, On the formation of singularities of solutions of
nonlinear differential systems in antistokes directions, Inv. Math.  45 (3) (2001) 425-485.

\bibitem{Duke} O. Costin, On Borel Summation and Stokes Phenomena for Rank-1 Nonlinear Systems of Ordinary
Differential Equations, Duke Math J. 93 (2) (1998) 289-344.

\bibitem{CT01} O. Costin, S. Tanveer, Existence and uniqueness for a class of nonlinear higher-order
partial differential equations in the complex plane, Comm. Pure Appl. Math. LIII (2000) 1092-1117.

\bibitem{CT02} O. Costin, S. Tanveer, Nonlinear evolution of PDEs in $\mathbb{R}^+ \times \mathbb{C}^2$:
existence and uniqueness of solutions, asymptotic and Borel summability properites, Annales De L'Institut
Henri Poicare' (C) Analyse Non Line'aire, 24 (5) (2007) 795-823.

\bibitem{CT03} O. Costin, S. Tanveer, Complex Singularity Analysis for a nonlinear PDE, Comm. PDEs, 31 (4) (2006) 593-637.

\bibitem{smalltime} O. Costin, S. Tanveer, Short Time Existence and Borel Summability in the Navier-Stokes Equation in $\mathbb{R}^3$, Communications in Paritial Differential Equations, 34 (8) (2009) 785-817. 

\bibitem{CT1} O. Costin, S. Tanveer, Nonlinear evolution PDEs in $\mathbb{R}^+ \times \mathbb{C}^d$: existence and uniqueness of solutions, asymptotic and Borel summability properties, Annales De L Institute Henri Poincare-Analyse Non Lineare, 24 (5) (2007) 795-823. 

\bibitem{scripta} O. Costin, G. Luo, S. Tanveer, An Integral Equation Approach To Smooth 3-D Navier-Stokes Solution, Physica Scritpa, T132 (014040) (2008)

\bibitem{longtime} O. Costin, G. Luo, S. Tanveer, Integral Formulation of 3-D Navier-Stokes and Longer Time Existence of Smooth Solutions, Comm. Contemp. Math, 13 (3) (2011) 407-462. 

\bibitem{Ecalle} J. Ecalle, Fonctions Resurgentes, Publications Mathematiques D' Orsay, (1981).

\bibitem{Ecalle1} J. Ecalle in Bifurcations nad periodic orbits of vector fields, NATO ASI Series, 408 (1993).

\bibitem{t-analyticity2} C. Foias, R. Temam, Some analytic and geometric properties of the solution of the evolution Navier-Stokes equation, J. Math. Pures Appl. (9) 58 (1979) 339-368. 

\bibitem{Boussinesq} T. Hou, C. Li, Global Well-Posedness of the Viscous Boussinesq Equations, Discrete and Continuous Dynamical Systems, 12 (1) (2005) 1-12.

\bibitem{t-analyticity1} G. Iooss, Application do la theorie des semi-groupes a l'etude de la stabilite des ecoule-ments laminaires, J. Mecanique, 8 (1969) 477-507.

\bibitem{Lutzetal} D.A. Lutz, M. Miyake, R. Schafke, On the 
Borel summability of divergent solutions of the heat equation, Nagoya
Math. J. 154 (1999) 1-29.

\bibitem{Raza} D. Razafindralandy, A. Hamdouni, Time Integration algorithm based on divergent
series summation, for ordinary and partial differential equations, J. Comp. Physics, 236 (2013) 56-73.

\bibitem{Thesis} H. Rosenblatt, Ph.D. Thesis, The Ohio State University, Asympototics and Borel Summability: Applications to MHD, Boussinesq Equations and Rigorous Stokes Constant Calculations (2013)

\bibitem{Temam} R. Temam, Navier-Stokes Equations: Theory and Numerical Analysis. AMS Chelsea Publishing, 2000.

\end{thebibliography}







\end{document}